\documentclass[12pt]{amsart}

\usepackage{amsfonts,amsmath,latexsym,amssymb,verbatim,amsbsy,times}
\usepackage{amsthm, esint}
\usepackage{graphicx}
\usepackage{caption}
\usepackage{subcaption}
\usepackage{dsfont}
\usepackage{tikz}
\usetikzlibrary{patterns}

\usepackage[margin=.75in]{geometry}
\usepackage[colorlinks=true, pdfstartview=FitV, linkcolor=blue,citecolor=green, urlcolor=blue]{hyperref}

\usepackage{enumitem}

\theoremstyle{plain}
\newtheorem{THEOREM}{Theorem}[section]
\newtheorem{COROL}[THEOREM]{Corollary}
\newtheorem{LEMMA}[THEOREM]{Lemma}
\newtheorem{PROP}[THEOREM]{Proposition}

\newtheorem{CONJECTURE}[THEOREM]{Conjecture}

\theoremstyle{definition}
\newtheorem{DEF}[THEOREM]{Definition}

\theoremstyle{remark}
\newtheorem{REMARK}{Remark}[section]

\newcommand{\N}{\ensuremath{\mathbb{N}}}   
\newcommand{\R}{\ensuremath{\mathbb{R}}}   

\def \one {{\mathds{1}}}

\def \a {\alpha}
\def \b {\beta}
\def \d {\delta}
\def \g {\gamma}
\def \e {\varepsilon}

\def \l {\lambda}

\def \s {\sigma}
\def \th {\theta}

\def \t {\tau}

\def \O {\Omega}

\def \cP {\mathcal{P}}

\def \wtA {\widetilde{A}}

\def \< {\langle}
\def \> {\rangle}
\def \p {\partial}

\DeclareMathOperator{\diam}{diam} %
\DeclareMathOperator{\supp}{supp} %

\def \cP {\mathcal{P}}

\def \dd  {\mathrm{d}}

\def \dt  {\, \mathrm{d}t}
\def \dt  {\, \mathrm{d}t}

\def \dy  {\, \mathrm{d}y}

\def \dr  {\, \mathrm{d}r}

\def \ds  {\, \mathrm{d}s}

\def \uphi  {\underline{\phi}}

\title[Cluster Formation in Alignment Dynamics]{Finite- and Infinite-Time Cluster Formation \\ for Alignment Dynamics on the Real Line}

\author{Trevor M. Leslie}
\address[T. Leslie]{Department of Applied Mathematics, Illinois Institute of Technology, Chicago, IL 60616}
\email{tleslie@iit.edu}

\author{Changhui Tan}
\address[C. Tan]{Department of Mathematics, University of South Carolina, Columbia, SC 29208}
\email{tan@math.sc.edu}



\usepackage{setspace}

\begin{document}

	\begin{abstract}
		We show that the locations where finite- and infinite-time clustering occurs for the 1D Euler-alignment system can be determined using only the initial data.  Our present work provides the first results on the structure of the finite-time singularity set and asymptotic clusters associated to a weak solution.  In many cases, the eventual size of the cluster can be read off directly from the flux associated to a scalar balance law formulation of the system.  
	\end{abstract}
	
\maketitle
\thispagestyle{empty} 


\setstretch{0.9}
\begin{small}
	\tableofcontents
\end{small}

\setstretch{1}

	\section{Introduction}
	
	\subsection{The Euler-alignment system} 
	This paper focuses on the 1-dimensional Euler-alignment system, which describes collective behavior among agents. It is a hydrodynamic analog of the celebrated Cucker--Smale system of ODE's, which we discuss below. The system is represented by the following equations:
	\begin{equation}
	\label{e:EA}
	\left\{ \begin{array}{rcl}
	\p_t \rho + \p_x(\rho u) & = & 0\,, \qquad (x,t)\in \R\times \R_+, \\
	\p_t (\rho u) + \p_x (\rho u^2) & = &  \displaystyle \int_{\R} \rho(x,t)\rho(y,t) \phi(x-y) (u(y,t) - u(x,t))\dy.
	\end{array}
	\right.
	\end{equation}
	The two functions $u$ and $\rho$ represent the velocity and nonnegative density profile of a group of agents. The communication protocol $\phi:\mathbb{R}\to \mathbb{R}$ describes the strength of the nonlocal alignment interactions. It is assumed to be nonnegative and even. The system \eqref{e:EA} is accompanied by initial data:
	\[
	\rho(x,0)=\rho^0(x),\quad u(x,0)=u^0(x).
	\]
	
The last decade has seen a rapid development in the theory for the well-posedness and asymptotic behavior of the Euler-alignment system \eqref{e:EA}.  As the nature of \eqref{e:EA} depends strongly on the behavior of the communication protocol $\phi$ near the origin, this theory necessarily breaks into several different cases.  When the communication protocol $\phi$ is \emph{strongly singular}, meaning that it is non-integrable at the origin, the alignment force on the right-hand side of \eqref{e:EA}$_2$ is known to exhibit dissipative properties and has a nonlinear regularizing effect. Consequently, the system takes on a parabolic character. Notably, research presented in \cite{do2018global,shvydkoy2017eulerian,shvydkoy2017eulerian2,shvydkoy2018eulerian} has demonstrated that the solution remains globally regular for any smooth initial data with $\inf \rho_0>0$.  In the case where the initial data contains vacuum, we refer to \cite{arnaiz2021singularity,tan2019singularity} for discussions on singularity formation and \cite{fabisiak2022inevitable} for existence of weak solutions.  Researchers have explored various extensions of the system, incorporating features such as pressure \cite{constantin2020entropy}, misalignment \cite{miao2021global}, and attraction-repulsion forces \cite{kiselev2018global}. While progress has been made in the context of the multi-dimensional system, it is less well-studied in comparison to the one-dimensional case. Some partial results have been presented in \cite{bai2022global,chen2021global,danchin2019regular,lear2021unidirectional,leslie2019weak,shvydkoy2019global}. 	

We are interested in the case where the communication protocol $\phi$ is less singular, specifically when it is integrable at the origin. In such instances, the alignment force exhibits a nonlocal damping effect, and the system adopts a hyperbolic character. A significant observation known as the \emph{critical threshold phenomenon} was first reported in \cite{tadmor2014critical}. This phenomenon highlights that the global regularity of solutions is contingent on the initial data: `subcritical' initial data lead to global well-posedness, while `supercritical' initial data result in finite-time singularity formations.
Numerous research works have been dedicated to determining the critical threshold conditions for the Euler-alignment system and related systems \cite{bhatnagar2021critical, carrillo2016critical, choi2019global, choi2019hydrodynamic, ha2018global, he2017global, leslie2020lagrangian, tan2020euler}.  
In the next subsection, we will conduct a brief survey on the findings specifically related to the 1-dimensional Euler-alignment system \eqref{e:EA}.

Another celebrated feature of the Euler-alignment system is its \emph{asymptotic flocking behavior}. This behavior emerges under the assumption that the communication kernel has a \emph{heavy tail}, indicated by the condition:
	\begin{equation}\label{eq:heavytail}
	\int_1^\infty\phi(r)\dr=\infty.
	\end{equation}
In such cases, solutions of the system converge to a flocking state as $t\to \infty$: the velocity $u$ aligns with its average value, while the density $\rho$ stabilizes into a traveling wave form:
\[
u(x,t)\to \bar{u}:=\frac{\int_\R (\rho^0 u^0)(x)\dd x}{\int_\R \rho^0\dd x},\qquad
\rho(x + \overline{u}t,t)\to\rho_{\infty}(x).
\]
Here $\rho_\infty$ is the asymptotic density profile, which carries important information about the emergent flocking phenomenon. However, the structure of $\rho_\infty$ is relatively less understood.
Relevant discussions on the flocking phenomenon can be found in works such as \cite{ha2009simple,leslie2019structure,leslie2023sticky,shvydkoy2017eulerian2,tadmor2014critical}.

	
	\subsection{Regular solutions and the critical threshold conditions}
	\label{ssec:CTC}
	
As mentioned above, the existence or non-existence of a global smooth solution $(\rho, u)$ to the 1-dimensional Euler-alignment system \eqref{e:EA} with smooth initial data $(\rho^0, u^0)$ can be determined from the \textit{critical threshold condition} (CTC). The conditions can be precisely characterized using an auxiliary quantity introduced in \cite{carrillo2016critical}:
\[
e^0(x) = \p_x u^0(x) + \phi*\rho^0(x), \qquad x\in \R.
\] 
Here `$*$' denotes convolution in the spatial variable.  Critical threshold conditions are available for two types of communication protocols:	
\begin{itemize}
	\item \emph{Bounded} communication.
	\item \emph{Weakly singular} communication: $\phi$ is unbounded but integrable at the origin. More precisely, we say $\phi$ is weakly singular with order $\beta\in(0,1)$ if there exist positive constants $R>0$ and $c>0$ such that the following lower bound holds:
	\begin{equation}\label{eq:phiws}
	\phi(r)\geq c\, r^{-\b},\quad\forall~r\in(0,R).
	\end{equation}
\end{itemize}

A sharp CTC has been established in \cite{carrillo2016critical} for bounded communication protocols, which says 
\begin{itemize}
	\item[I.] If $e^0(x)\geq0$ for all $x\in\R$, then the solution is globally regular.
	\item[II.] If $e^0(x_0)<0$ for some $x_0\in\R$, then the solution develops a singularity at $x(x_0,T_*)$ for some finite time $T_*$.
\end{itemize}
In this context, we use the notation $x(a,t)$  for the characteristic path originating from $a$, which satisfies the following ordinary differential equation:
\[
\dot{x}(a,t) = u(x(a,t),t),\quad x(a,0) = a.
\]
For bounded protocols, these conditions are often referred to as `subcritical' (I) and `supercritical' (II).

When the communication protocol is weakly singular, another type of finite time singularity formation was discovered in \cite{tan2020euler}. Under the subcritical CTC (I), we have the following.
\begin{itemize}
	\item[III.] If $e^0(x_0)=0$ for some $x_0\in\R$, then the behavior depends on the type of communication protocol:
	\begin{itemize}
	\item[(i)] If $\phi$ is bounded, then the solution is globally regular.
	\item[(ii)] If $\phi$ is weakly singular, then the solution might develop a singularity at $x(x_0,T_*)$ for some finite time $T_*$. 	
	\end{itemize}		
\end{itemize}
An example was provided in \cite{tan2020euler} to illustrate the singularity formation, assuming $\phi$ has a heavy tail.
Note that this is the only scenario where the two types of protocols can lead to different behaviors. We may refer it as the critical case.

The nature of the singularity in II and III(ii) is commonly known as a \emph{singular shock}, arising when two characteristic paths collide and lead to a shock discontinuity in velocity. Additionally, there is a concentration of mass, which we refer to as \emph{clustering} in this paper.
	
In \cite{leslie2020lagrangian}, the author proved a refinement of the above results that renders the CTC more meaningful in the presence of vacuum.
When $x\notin \supp \rho^0$, the physical velocity $u^0(x)$ is undefined, making $e^0(x)$ ill-defined. To address this, the author introduced an anti-derivative of $e^0$, denoted as:
\begin{equation}
	\psi^0 = u^0 + \Phi*\rho^0, \qquad \Phi(x) = \int_0^x \phi(r)\dr.
\end{equation}
The CTC are then expressed in terms of the monotonicity of $\psi^0$ inside the support of $\rho^0$. For instance, if $\phi$ is bounded, the solution is globally regular if and only if $\psi^0$ is nondecreasing in $\supp(\rho^0)$.  The theory presented in \cite{leslie2020lagrangian} was further developed in \cite{lear2022geometric} through a comprehensive study of the evolution of characteristic paths; the latter will inform the following discussion.

The characteristics are simplest for the degenerate protocol $\phi \equiv 0$, for which the Euler-alignment system \eqref{e:EA} reduces to the well-studied pressureless Euler system.  For pressureless Euler, the characteristic paths of classical solutions are always straight lines, leading to three possible scenarios: (I) separation linearly in time, (II) collision in finite time, or (III) running parallel for all time. However, the introduction of the alignment force leads to intriguing new asymptotic behaviors.

For bounded and heavy-tailed communication protocols, assuming $\psi^0$ is nondecreasing in $\text{supp}(\rho^0)$, it was demonstrated in \cite{lear2022geometric} that for any $a, b\in\text{supp}(\rho^0)$ with $a<b$, the characteristic paths $x(a,t)$ and $x(b,t)$ are globally well-defined, and their distance satisfies the following quantitative bounds (for some positive constants $c$ and $C$ that do not depend on $a$, $b$, or $t$):
\begin{equation}
	\label{e:xsimpsi}
	c(\psi^0(b) - \psi^0(a))
	\le \lim_{t\to \infty} \big( x(b,t) - x(a,t) \big) 
	\le C(\psi^0(b) - \psi^0(a)).
\end{equation}
The comparison \eqref{e:xsimpsi} was then utilized to study the structure of the asymptotic density profile $\rho_\infty$. One especially interesting situation occurs when $\psi^0(a) = \psi^0(b)$, in which case the distance between characteristic paths tends to zero as $t\to +\infty$:
\[\lim_{t\to \infty} \big( x(b,t) - x(a,t) \big)=0.\]
If there is any mass trapped between the two converging characteristics, a concentration of mass will develop in the asymptotic density profile $\rho_\infty$. This phenomenon is referred to as \emph{infinite-time clustering}, and it represents a distinctive feature of the alignment interaction. A discrete version of this phenomenon was previously investigated in \cite{ha2019complete} for the classical Cucker--Smale system (which we note often exhibits substantially different behavior than its sticky particle version considered below, c.f. Remark \ref{r:1.1}). A particular bi-cluster formation was studied in \cite{cho2016emergence}.

We summarize the key results cited above in simplified form as follows:
\begin{PROP}\label{p:classical}
Suppose $\rho^0 \in C(\R)$ is compactly supported and $u^0\in C^1(\R)$.  Let $\phi$ be nonnegative, even, and locally integrable.  Assume $a<b$.  
\begin{itemize}
	\item [\textnormal{I.}] Suppose $\psi^0$ is nondecreasing and $\psi^0(a)<\psi^0(b)$. Then there exists $c>0$ such that 
			\[
			x(b,t) - x(a,t)\ge c>0 
			\qquad \text{ for all } t\ge 0.
			\]
			\item [\textnormal{II.}] Suppose $\psi^0(a)>\psi^0(b)$. Then $u$ loses regularity in finite time.
			\item [\textnormal{III.}] Suppose $\psi^0$ is nondecreasing and $\psi^0(a) = \psi^0(b)$.  
			\begin{itemize}
				\item [\textnormal{(i)}] If $\phi$ is bounded, then $x(b,t)>x(a,t)$ for all $t\in [0,\infty)$.
				\item [\textnormal{(ii)}] If $\phi$ is bounded and heavy-tailed, then we have we have $x(b,t) - x(a,t)\to 0$ as $t\to \infty$.
				\item [\textnormal{(iii)}] If $\phi$ is weakly singular and heavy-tailed, and if $\rho^0\big|_{[a,b]}\not\equiv 0$, then finite-time blowup occurs.
			\end{itemize}
		\end{itemize} 	
	\end{PROP}
	
The primary objective of this study is to investigate the phenomena of finite- and infinite-time clustering and to derive predictions from the initial data. We seek to extend the current theory by addressing the following gaps in Proposition \ref{p:classical}:
\begin{itemize}
	\item For finite-time clustering: The current results only demonstrate finite-time singularity formation in cases II and III(iii). However, the behavior after the clusters form remains unknown as the classical solution ceases to exist. We aim to understand the subsequent evolution of these clusters once they have formed.
	\item For infinite-time clustering: The existing results treat the infinite-time clustering phenomenon only when the solution is globally regular, meaning that there is no finite-time clustering. Moreover, the simple characterization of the infinite-time clusters in terms of the monotonicity properties of $\psi^0$ is not expected to survive in the presence of finite-time clusters; we seek a suitable generalization of this characterization.
\end{itemize}

To achieve our objective, we need to consider an appropriate class of weak solutions. For this purpose, we will build upon the theory recently established in \cite{leslie2023sticky}, which provides a valuable framework for our study. It will allow us to explore the intricate dynamics of the Euler-alignment system and investigate the finite- and infinite-time clustering phenomena in a comprehensive manner. 
	
	

\subsection{Weak solutions and the scalar balance law}
The well-posedness theory for weak solutions to systems of conservation laws poses significant challenges, particularly concerning uniqueness. In \cite{leslie2023sticky}, the authors establish the well-posedness theory for weak solutions to \eqref{e:EA} by employing the approach of Brenier and Grenier \cite{brenier1998sticky} on the 1D pressureless Euler system ($\phi\equiv0$) to address this issue. The key idea is to reduce the system to a single scalar balance law:
\begin{equation}\label{e:SB}
	\p_tM + \p_x(A(M)) = (\Phi*\p_x M)\p_x M, \qquad \Phi(x) = \int_0^x \phi(r)\dr,
\end{equation}
supplemented with initial conditions $M(\cdot, 0) = M^0$.  Here $M(t):\R\to [-\frac12, \frac12]$ is the cumulative distribution function for the density $\rho(t)$ (shifted by a constant for technical reasons) and the flux $A$ is determined from $\rho^0$ and $u^0$.  Unlike the situation for \eqref{e:EA}, it is fairly straightforward to establish (though not trivial to justify) entropy conditions for \eqref{e:SB} that are sufficient to guarantee uniqueness (c.f. \cite{leslie2023sticky}).  We therefore refer to the solution $(\rho, u)$ of \eqref{e:EA} that we recover from \eqref{e:SB} as the \textit{entropy solution} of \eqref{e:EA}.  
	
	Let us give an extremely brief description of the construction of an entropy solution to \eqref{e:EA}.
	Starting with initial data $(\rho^0, u^0)\in \cP_c(\R)\times L^\infty(\dd\rho^0)$ (where $\cP_c(\R)$ denotes the space of compactly supported probability measures on $\R$), we define the corresponding cumulative  distribution function $M^0$ and its generalized inverse $X^0$ as follows
	\begin{equation}\label{eq:M0}
	M^0(x) = -\frac12 + \rho^0((-\infty, x]),\quad X^0(m) = \inf\left\{x:M^0(x)\ge m\right\}.
	\end{equation}
	We define the flux $A:[-\frac12, \frac12]\to \R$ of the scalar balance law \eqref{e:SB} as
	\begin{equation}
	\label{e:defA}
	A(m) = \int_{-\frac12}^m \psi^0 \circ X^0(\widetilde{m})\,\dd\widetilde{m},\quad \psi^0=u^0+\Phi\ast \rho^0.
	\end{equation}
	Note that $A$ is Lipschitz, with $A(-\frac12) = 0$.  Having determined $(M^0, A)$ from $(\rho^0, u^0)$, we evolve the scalar balance law \eqref{e:SB}.  There is a unique entropy solution, which we denote by $M = M(x,t)$,  associated to the initial data $M^0$ and flux $A$.  We generate this entropy solution through a front-tracking approximation scheme; remarkably, the positions and magnitudes of our fronts can be encoded using the Cucker--Smale dynamics, supplemented with completely inelastic collision rules.  (We will describe this \textit{sticky particle Cucker--Smale} system in detail in Section \ref{sec:preliminaries}, and it will play a crucial role in the proof of our main theorem.) Finally, we recover the solution to \eqref{e:EA} via $\rho=\p_xM$ and $P=-\p_tM$. The velocity $u$ is then the Radon-Nikodym derivative of the measure $P$ with respect to $\rho$. It is shown in \cite{leslie2023sticky} that this pair $(\rho, u)$ solves the Euler-alignment system \eqref{e:EA} in the sense of distributions and satisfies the initial data in an appropriate sense.  
	
	We also set notation for the generalized inverse $X(t)$ of $M(t)$, which we will use extensively below.
	\begin{equation}\label{eq:X}
	X(m,t) = \inf\left\{x\in \R: M(x,t)\ge m\right\}.
	\end{equation}
	Note that $M^0$ and $M(t)$ are nondecreasing, right-continuous, and defined on $\R$, while $X^0$ and $X(t)$ are nondecreasing, left-continuous, and defined on $(-\frac12, \frac12]$.  We refer to elements of $(-\frac12, \frac12]$ as \emph{mass labels}.

	\subsection{Terminology and notation for cluster formation}
	The focus of this paper is on cluster formation of mass labels in solutions of the Euler-alignment system \eqref{e:EA}. In particular, we would like to describe the phenomena of finite- and infinite-time clustering by examining the initial conditions.  
	
	We begin by defining our clusters in terms of the function $X(\cdot, t)$, which encodes the `location' of each mass label. 
	
	\begin{DEF}[Clusters]
		\label{def:cluster}
		Given $m\in (-\frac12, \frac12]$, we say that there is a \emph{$t$-cluster at $m$} (or a \emph{finite-time cluster at $m$}, if the value of $t$ is not important) if there exists $m'<m$ such that $X(m',t) = X(m,t)$.  In this case, the $t$-cluster at $m$ is defined to be the largest interval of the form $(m',m'']$, containing $m$, such that $X(\cdot, t)$ is constant on $(m',m'']$.  We refer to $0$-clusters as \textit{initial clusters}.
		
		We say there is an \emph{infinite-time cluster at $m$} if there exists $m'<m$ such that $X(m,t) - X(m',t)\to 0$ as $t\to \infty$.  In this case, the infinite-time cluster at $m$ is defined to be the largest interval of the form $I = (m',m'']$ or $I = (m',m'')$, containing $m$, such that $\diam X(I,t)\to 0$ as $t\to \infty$. 
	\end{DEF}
	
	It is worth noting that if $(m',m'']$ is the $t$-cluster at $m$, then $X(\cdot, t)$ is constant on $(m',m'']$ by Definition~\ref{def:cluster}, but either of the possibilities $X(m',t) = X(m'', t)$ or $X(m', t)\ne X(m'',t)$ may occur. Nevertheless, we feel that the left-continuity of $X(\cdot, t)$ provides us with a compelling reason to define our finite-time clusters as half-open intervals.
		
	Our definition above, and the analysis of \cite{leslie2023sticky}, guarantee that existing clusters cannot `unstick.' Hence, the size of a cluster can only grow in time, and if there is a finite-time cluster at $m$, then there is an infinite-time cluster at $m$ as well, with the infinite-time cluster containing the finite-time cluster.

  We say that that there is \emph{no finite-time clustering at $m$} if there is no $t$-cluster at $m$ for all $t\in [0,\infty)$, and we say that there is \textit{no infinite-time clustering at $m$} if there is no infinite-time cluster at $m$.
	
	In what follows, we use $A^{**}$ to denote the lower convex envelope of the flux $A$.  We will give a brief reminder of some definitions related to convexity in Section \ref{ssec:convexity}.  
	
	\begin{DEF}\label{def:Sigma} We break up the interval $(-\frac12, \frac12]$ into disjoint regions, on which we will observe different clustering behavior, as follows:
		\begin{itemize}[itemsep = .5em]
			\item We define the \textit{subcritical region} $\Sigma_+$ by
			\[
			\begin{split} 
			\Sigma_+ 
			& = \{m\in (-\tfrac12, \tfrac12]: A^{**} \text{ is not linear on any interval of the form } (m',m]\}.
			\end{split} 
			\]
			\item We define the \textit{critical region} $\Sigma_0$ by 
			\[
			\Sigma_0 = \bigcup_{\substack{ A \text{ is linear and } \\ \text{ equal to } A^{**} \text{ on } (m',m'']}} (m',m''].	
			\]
			\item We define the \textit{supercritical region} $\Sigma_-$ as the (open) set on which $A>A^{**}$:
			\[
			\Sigma_- = \{m\in (\tfrac12, \tfrac12): A(m)>A^{**}(m)\}.
			\]  	
		\end{itemize}
	\end{DEF}	
	
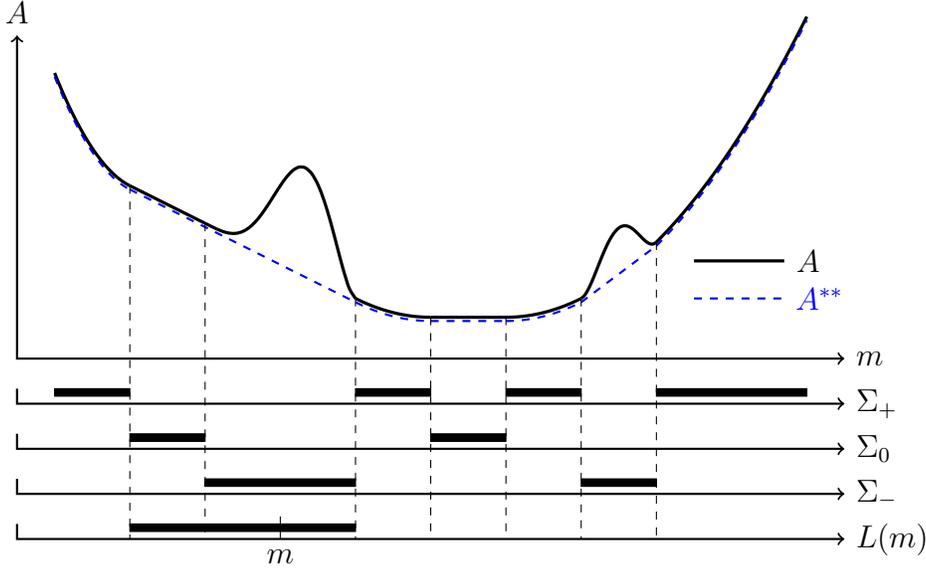
\begin{figure}[ht]
		\begin{tikzpicture}
		\draw[thick,<->] (10.5,.2) node[right]{$m$} -- (-.5,.2) -- (-.5,4.5) node[above]{$A$};
		
		\draw[domain=0:1, smooth, variable=\x, very thick] plot ({\x}, {\x^2-5/2*\x+4});
		\draw[very thick] (1, 2.5) -- (2,2);
		\draw[domain=2:4, smooth, variable=\x, very thick] plot ({\x}, {(\x-2)^4*(\x-4)^2+3-\x/2});
		\draw[domain=4:5, smooth, variable=\x, very thick] plot ({\x}, {(5-\x)^2/4+3/4});
		\draw[very thick] (5, .75) -- (6,.75);
		\draw[domain=6:7, smooth, variable=\x, very thick] plot ({\x}, {(\x-6)^2/4+3/4});
		\draw[domain=7:8, smooth, variable=\x, very thick] plot ({\x}, {(\x-6)^2/4+3/4+10*(\x-7)^2*(\x-8)^2});
		\draw[domain=8:10, smooth, variable=\x, very thick] plot ({\x}, {(\x-6)^2/4+3/4});
		
		\draw[domain=0:1, smooth, variable=\x, thick, dashed, blue] plot ({\x}, {\x^2-5/2*\x+4-.05});
		\draw[thick, dashed, blue] (1, 2.45) -- (4,.95);
		\draw[domain=4:5, smooth, variable=\x, thick, dashed, blue] plot ({\x}, {(5-\x)^2/4+3/4-.05});
		\draw[thick, dashed, blue] (5, .7) -- (6,.7);
		\draw[domain=6:7, smooth, variable=\x, thick, dashed, blue] plot ({\x}, {(\x-6)^2/4+3/4-.05});
		\draw[thick, dashed, blue] (7, .95) -- (8,1.7);
		\draw[domain=8:10, smooth, variable=\x, thick, dashed, blue] plot ({\x}, {(\x-6)^2/4+3/4-.05});
		
		\draw[thick,->] (-.5,-.2) -- (-.5,-.4) -- (10.5,-.4) node[right]{$\Sigma_+$};
		\draw[fill] (0,-.3) rectangle (1,-.2);
		\draw[fill] (4,-.3) rectangle (5,-.2);
		\draw[fill] (6,-.3) rectangle (7,-.2);
		\draw[fill] (8,-.3) rectangle (10,-.2);
		
		\draw[thick,->] (-.5,-1.4) -- (-.5,-1.6) -- (10.5,-1.6) node[right]{$\Sigma_-$};
		\draw[fill] (2,-1.5) rectangle (4,-1.4);
		\draw[fill] (7,-1.5) rectangle (8,-1.4);
		
		\draw[thick,->] (-.5,-.8) -- (-.5,-1) -- (10.5,-1) node[right]{$\Sigma_0$};
		\draw[fill] (1,-.9) rectangle (2,-.8);
		\draw[fill] (5,-.9) rectangle (6,-.8);

		\draw[thick,->] (-.5,-2) -- (-.5,-2.2) -- (10.5,-2.2) node[right]{$L(m)$};
		\draw[fill] (1,-2.1) rectangle (4,-2);
		\draw (3,-1.9) -- (3,-2.2) node[below]{$m$};

		\draw[dashed] (1,2.5) -- (1,-2.2) (2,2) -- (2,-2.2) (4,1) -- (4,-2.2) (5,.75) -- (5,-2.2) (6,.75) -- (6,-2.2) (7,1) -- (7,-2.2) (8,1.75) -- (8,-2.2);
		
		\draw[very thick] (8.5,1.5) -- (9.7,1.5) node[right]{$A$};
		\draw[thick,blue,dashed] (8.5,1) -- (9.7,1) node[right]{$A^{**}$};
		\end{tikzpicture}
		\caption{An illustration of the three regions of Definition \ref{def:Sigma}, and an example of a set $L(m)$ that contains multiple connected intervals in $\Sigma_0$ and $\Sigma_-$. 
		}\label{fig:AA**}
	\end{figure}

	Let us provide some brief commentary on the definition above. Neglecting a set of measure zero, we can understand the three regions in terms of the table below, which provides a more intuitive picture.  (Note that primes on $A$ and $A^{**}$ will always represent derivatives below.) We stress, however, that the formulation of Definition \ref{def:Sigma} is better suited for our analysis, as will become clear below.
	\begin{center}
	\def\arraystretch{1.5}
	\begin{tabular}{|c|c|c|}
	\hline
	&$A(m)=A^{**}(m)$&$A(m)>A^{**}(m)$\\ \hline
	$(A^{**})''(m)>0$&$\Sigma_+$&$\emptyset$\\ \hline
	$(A^{**})''(m)=0$&$\Sigma_0$&$\Sigma_-$\\ \hline	\end{tabular}
	\end{center}
	See Figure \ref{fig:AA**} for an illustration of the three regions. 
	
	Let us also note that in the `typical' case where each of $\Sigma_-$, $\Sigma_0$, and $\Sigma_+$ has finitely many connected components, both $\Sigma_+$ and $\Sigma_0$ are unions of half-open intervals. Consequently, the disjoint union $\Sigma_+ \cup \Sigma_0 \cup \Sigma_-$ contains all points in $(-\frac{1}{2}, \frac{1}{2}]$ \textit{except} for the right endpoints of the intervals that constitute the connected components of $\Sigma_-$. This `missing' set necessarily has Lebesgue measure zero (and one of our assumptions below will actually force it to be finite for the cases we consider).		
	
	We set one more piece of notation before moving on.
	\begin{DEF}
		\label{def:Lm}
		For $m\notin \Sigma_+$, we define the set $L(m)$ to be the largest half-open interval $(m',m'']$ containing $m$, such that $A^{**}$ is linear on $L(m)$.  For $m\in \Sigma_+$, we define $L(m) = \{m\}$.	
	\end{DEF}
	
	Ignoring endpoints, we may simply view $L(m)$ (for $m\notin \Sigma_+$) as the largest interval containing $m$ on which $(A^{**})''~=~0$. It may contain connected components of both $\Sigma_0$ and $\Sigma_-$, as illustrated in Figure \ref{fig:AA**}. 
		
		Let us briefly discuss how some of the results of Section \ref{ssec:CTC} translate to the present framework, and (at an informal level) how we plan to extend them. For simplicity, we assume in this discussion that our initial data is sufficiently regular, i.e., $\psi^0$ and $X^0$ are differentiable and $X^0$ is strictly increasing.  In this context, we have
		\[A''(m)=\frac{\dd}{\dd m}\psi^0(X^0(m))=e^0(X^0(m))\cdot(X^0)'(m).\]
		Then, we have following equivalent representation of CTC:
		\begin{align*}
			\text{I.}&\quad e^0(x)\geq0 \text{ for all } x\in\R  \iff A''(m)\geq0 \text{ for all } m\in(-\tfrac12,\tfrac12] \iff \Sigma_-=\emptyset,\\
			\text{II.}&\quad e^0(x_0)<0 \text{ for some } x_0\in\R  \iff A''(m_0)<0 \text{ for some } m_0\in(-\tfrac12,\tfrac12] \iff \Sigma_-\neq\emptyset.
		\end{align*}
		Furthermore, under the subcritical CTC (I), or equivalently, $\psi^0$ being nondecreasing, for $m'<m''$, we have the following equivalences for the assumptions in Proposition \ref{p:classical}:
		\begin{align*}
			\text{I.}\quad \psi^0(X^0(m')) < \psi^0(X^0(m'')) & \iff A'(m')<A'(m'') \iff A \text{ is not linear on } [m',m'']\\
			& \iff L(m')\neq L(m''),\\
			\text{III.}\quad \psi^0(X^0(m')) = \psi^0(X^0(m'')) & \iff A'(m')=A'(m'') \iff A \text{ is linear on } [m',m'']\\
			& \iff L(m') = L(m'').
		\end{align*}		
		Comparing the above discussion with Proposition \ref{p:classical}, we arrive at the following conjectures:
		\begin{CONJECTURE}\label{conj:comp}
			Mass labels from distinct $L(m)$'s never belong to the same cluster.  
		\end{CONJECTURE}
 \begin{CONJECTURE}\label{conj:comp2}
		If $\phi$ is heavy-tailed, then mass labels from the same $L(m)$ belong to the same infinite-time cluster.  
	\end{CONJECTURE}
		We will demonstrate later that these conjectures are indeed correct. More remarkably, we show they hold even when $\Sigma_-\neq\emptyset$. Indeed, the equivalences above hold without the assumption $\Sigma_-=\emptyset$, if we replace $A$ by $A^{**}$ and we assume that $m',m''\notin \Sigma_-$.
		
	Let us also provide commentary on part II of Proposition \ref{p:classical}. The condition $\psi^0(a)>\psi^0(b)$, $a<b$   implies that finite-time cluster formation occurs `somewhere between $a$ and $b$.' However, the monotonicity of $\psi^0$ itself does not provide complete information about the clusters and their evolution over time. (For instance, $\psi^0(X^0(m'))>\psi^0(X^0(m''))$, $m'<m''$ implies nothing about whether $m'$ and $m''$ eventually belong to the same cluster.) Our framework allows us to obtain more detailed information about finite-time clustering. One new finding is:
		\begin{center} 
			\emph{Any two mass labels from the same connected component of $\Sigma_-$ cluster together in finite time}.  
		\end{center}
	We present precise versions of these statements in Theorem \ref{t:regions} below, along with other features that extend Proposition \ref{p:classical}.

	\subsection{Global Assumptions and Summary of Results}
	
	Throughout our manuscript, we fix initial data $(\rho^0, u^0)$ and make use of the following global assumptions, except where explicitly stated otherwise:
	\begin{itemize}
		\item [(A1)] The communication protocol $\phi$ is locally integrable, even, and radially nonincreasing.
		\item [(A2)] $(\rho^0, u^0)\in \cP_c(\R)\times L^\infty(\dd\rho^0)$, and $(\rho, u)$ is the associated entropy solution.
		\item [(A3)] $M^0$, $X^0$, $A$, and $X$ are defined as in \eqref{eq:M0}, \eqref{e:defA}, \eqref{eq:X}, and the sets $\Sigma_-$, $\Sigma_0$, $\Sigma_+$, and $L(m)$ are defined accordingly as in the previous subsection.
		\item [(A4)] $A$ is convex in a neighborhood of every point $m$ of the boundary $S:=\p\Sigma_-$ of $\Sigma_-$.
	\end{itemize}
	The technical assumption (A4) is the only one we have not previously mentioned. It will simplify the structure of the set $\Sigma_-$, as will be spelled out below in Lemma \ref{lem:S}; it will also play a key role in our analysis of the clustering behavior in $\Sigma_-$.  
	
	We now state our main result.  
	\begin{THEOREM}
		\label{t:regions}
		Fix $m\in (-\frac12, \frac12]$.  The following statements describe the clustering behavior at $m$. 
		\begin{itemize}
			\item [\textnormal{I.}] If $m\in \Sigma_+$, then there is no finite- or infinite-time clustering at $m$.
			\item [\textnormal{II.}] If $m\in \Sigma_-$, then there is a $t$-cluster at $m$ for all sufficiently large $t\ge 0$.  Moreover, if $(m_-, m_+)$ is the connected component of $\Sigma_-$ containing $m$, then for any $\widetilde{m}\in (m_-, m_+)$ there exists a time $T \ge 0$ such that $m$ and $\widetilde{m}$ lie in the same $t$-cluster for all $t\ge T$.  
			\item [\textnormal{III.}]  
			\begin{itemize}
				\item [\textnormal{(i)}] Suppose $\phi$ is bounded.  If $(m',m'']$ is a finite-time cluster at $m$, then either $m\in \Sigma_0$ and $(m',m'']$ is an initial cluster, or $m\in \Sigma_-$ and $(m',m'']\subseteq (m_-, m_+]$, where $(m_-, m_+)$ is the connected component of $\Sigma_-$ containing $m$.  No other finite-time clusters are possible.
				\item [\textnormal{(ii)}] Suppose $\phi$ is heavy-tailed.  If $m\notin \Sigma_+$, then there is an infinite-time cluster at $m$, and it is equal to $L(m)$.    
				\item [\textnormal{(iii)}] Suppose $\phi$ is heavy-tailed and weakly singular.  If $m\notin \Sigma_+$, then there exists a finite time $T$ such that $L(m)$ is a $t$-cluster at $m$ for all $t\ge T$.  
			\end{itemize}
		\end{itemize}
	\end{THEOREM}

\begin{REMARK}
	\label{r:1.1}
The theorem above is clearly a significant upgrade over Proposition \ref{p:classical}, which summarized the most relevant results on classical solutions.  Indeed, our new theorem satisfactorily addresses the gaps outlined in Section \ref{ssec:CTC} and applies to very general initial conditions.  Let us clarify how it relates to two other lines of research that we have not yet emphasized.

We briefly mentioned above the paper \cite{ha2019complete}, where the authors gave a comprehensive study of the limiting configurations associated to the 1D Cucker--Smale system with bounded communication.  Their analysis is concerned with the `free-flow' dynamics, where agents follow the Cucker--Smale ODE's for all time, without modification for the occurrence of collisions.  Their analysis, comprehensive though it is for the classical Cucker--Smale system, has no hydrodynamic analog except for the case when no collisions occur.  (The previously discussed work \cite{leslie2019structure} is in this direction.)  Once collisions occur, the proper way to pass to a hydrodynamic limit is through the \textit{sticky particle} Cucker--Smale dynamics, as the present authors showed in our previous paper \cite{leslie2023sticky}.  As soon as collisions are allowed, the sticky particle Cucker--Smale dynamics can differ drastically from those of their classical counterpart.

The one case where the free-flow dynamics \textit{do} have substantial bearing on the sticky particle dynamics is in the setting of degenerate communication $\phi \equiv 0$, where \eqref{e:EA} reduces to the 1D pressureless Euler equations.  In this degenerate case, Theorem \ref{t:regions} can be recovered from already-existing theory; see for example \cite{huangwang2001, natile2009wasserstein, brenier2013sticky, cavalletti2015simple}. (The paper \cite{brenier2013sticky} considers a more general system which includes pressureless Euler and Euler--Poisson as special cases.)  These works all rely implicitly or explicitly on the fact that the `sticky particle' dynamics in that setting can be recovered from the `free-flow' dynamics using a certain $L^2$ projection onto the convex cone of nondecreasing functions.  (This is essentially the key observation of \cite{natile2009wasserstein}.)  However, simple counterexamples show that the sticky	particle Cucker--Smale dynamics, discussed below, cannot be recovered as such a projection.  Consequently, the techniques of the above-cited works do not appear to apply to our problem.  In particular, even at the discrete level, the result of applying the $L^2$ projection to the analysis of Ha et.$  $ al \cite{ha2019complete} is in general unrelated to the cluster formation described in Theorem \ref{t:regions}.
\end{REMARK}

	\subsection{Outline of the paper}	
	Our proof relies on a discretization procedure from \cite{leslie2023sticky} involving the so-called `sticky particle Cucker--Smale' dynamics; in Section \ref{sec:preliminaries}, we review the salient parts of this procedure for the convenience of the reader.  We also establish conventions and review some standard facts related to convex functions.  The remaining sections \ref{sec:sub}, \ref{sec:super}, and \ref{sec:phidependent} contain all the new analysis.  These sections essentially track the statements I, II, and III, respectively, of Theorem \ref{t:regions}, with some caveats spelled out below.  We give a slightly more detailed summary presently.
	
	In Section \ref{sec:sub}, we prove Conjecture \ref{conj:comp}: mass labels from distinct $L(m)$'s cannot belong to the same finite- or infinite-time cluster.  Since $L(m)$ is a singleton for every $m\in \Sigma_+$, this implies statement I of Theorem \ref{t:regions} as a special case.  The greater generality of this framework will pay dividends in the proofs of statement II and (especially) statement III.  
	
	Section \ref{sec:super} concerns the supercritical region $\Sigma_-$.  We give the proof of part II of Theorem \ref{t:regions}, which relies strongly on assumption (A4).  We also prove that if $m$ lies inside the connected component $(m_-, m_+)$ of $\Sigma_-$, then no $t$-cluster at $m$ can extend beyond $(m_-, m_+]$ unless it contains all of $(m_-, m_+]$.   Roughly speaking, this tells us that the interval $(m_-, m_+]$ can be treated as an indivisible unit for purposes of the larger-scale clustering analysis.  
	
	In Section \ref{sec:phidependent}, we prove part III of Theorem \ref{t:regions} (including Conjecture \ref{conj:comp2}, and more).  This is the only place in the paper where we specialize our assumptions on $\phi$ beyond (A1).  The previously mentioned statement on the `indivisibility' of the connected components of $\Sigma_-$ becomes extremely useful for proving III(i); in particular it is essential when dealing with the situation where a connected component of $\Sigma_-$ borders an initial cluster or another connected component of $\Sigma_-$.  Finally, the framework we have developed by the time we reach the proofs of statements III(ii) and III(iii) allows us to write the latter as adaptations of known arguments.

	\section{Preliminaries}
	
	\label{sec:preliminaries}
	
	In this section, we collect some preliminary results on the 1D Euler-alignment system \eqref{e:EA} and the corresponding role of the `sticky particle Cucker--Smale' dynamics. 
	
	First, however, we recall a few standard definitions and facts related to convex functions of a single variable, and we set a few conventions. We also justify our previous claim that $\Sigma_-$ has only finitely many connected components, as a consequence of the convexity assumption (A4).
	
	\subsection{Convex functions}
	\label{ssec:convexity}
	\begin{DEF}[Convexity]\label{def:conv} 
		Let $\O$ be an interval in $\R$.  We say that $A:\O\to \R$ is \textit{convex} if for every $m',m''\in \O$ and every $\l\in (0,1)$, we have 
		\[
		A((1-\l)m' + \l m'')\le (1-\l)A(m') + \l A(m'')
		\]
		We say $A$ is \textit{strictly convex} if the above inequality is  strict for all $m',m''\in \O$ and all $\l\in (0,1)$.  If $I\subseteq \O$ is a subinterval of $\O$, then we say that $A$ is (strictly) convex on $I$ if the restriction of $A$ to $I$ is (strictly) convex.  Finally, we say that $A$ is (strictly) convex in a neighborhood of $m\in \O$ if there exists an interval $I$, open relative to $\O$, such that $m\in I$ and $A$ is (strictly) convex on $I$.  
		
		The \textit{lower convex envelope} of $A:\O\to \R$ is the largest convex function $A^{**}:\O\to \R$ whose graph lies below that of $A$:
		\[ A^{**}(m) = \sup\left\{\wtA(m)~|~ \wtA \text{ is convex and } \wtA\leq A \text{ on all of } [-\tfrac12,\tfrac12]\right\}.\]
	\end{DEF}
	
	We also recall the following very useful Lemma, which is a direct consequence of the definition. 
	\begin{LEMMA}
		\label{l:convsec}
		Suppose $A$ is convex in an interval $I$. Let $[\widetilde{m}',\widetilde{m}'']$ and $[m',m'']$ be two sub-intervals of $I$ such that $\widetilde{m}'\leq m'$ and $\widetilde{m}''\leq m''$. Then
		\begin{equation}\label{eq:convsec}
		\frac{A(\widetilde{m}'')-A(\widetilde{m}')}{\widetilde{m}''-\widetilde{m}'}
		\le \frac{A(m'')-A(m')}{m''-m'}.	
		\end{equation}
	\end{LEMMA}

	The following elementary lemma will be used several times in Section \ref{sec:sub}.
	
	\begin{LEMMA}
		\label{l:A**segments}
		Let $A$ be a real-valued function defined on an interval containing $[m',m'']$.  Suppose that for some $\th\in (m',m'')$, it holds that
		\begin{equation}
		\label{e:thxi}
		(m',A(m')), \; (\th, A^{**}(\th)),\; \text{ and } (m'', A(m'')) \text{ are collinear.}
		\end{equation}
		Then $A(m') = A^{**}(m')$ and $A(m'') = A^{**}(m'')$, and consequently, $A^{**}$ is linear on $[m',m'']$.
	\end{LEMMA}
	\begin{proof}
		Pick $\l\in (0,1)$ such that $\th = (1-\l)m' + \l m''$.  Then 
		\begin{align*}
		(1-\l)A(m') + \l A(m'') & = A^{**}(\th) \\
		& \le (1-\l) A^{**}(m') + \l A^{**}(m'') \\
		& \le (1-\l) A(m') + \l A(m'').
		\end{align*}	
		Since the left and right sides are equal, this forces 
		\[
		(1-\l)(A(m') - A^{**}(m')) + \l (A(m'') - A^{**}(m'')) = 0.
		\]
		Since $\l\in (0,1)$ and $A-A^{**}\ge 0$, we must have $A(m') = A^{**}(m')$ and $A(m'') = A^{**}(m'')$.  Linearity of $A^{**}$ on $[m',m'']$ then follows from Lemma \ref{l:convsec}.
	\end{proof}

	The convexity assumption (A4) guarantees that the supercritical region~$\Sigma_-$ has only finitely many connected components, as we establish presently.
	
	\begin{LEMMA}\label{lem:S}
		The set $\Sigma_-$ has only finitely many connected components.
	\end{LEMMA}
	\begin{proof}
		Suppose not; then the set $S_r$ consisting of the right endpoints of the connected components of $\Sigma_-$ must have a limit point $m$ in $[-\frac12, \frac12]$.  Let $(m_j)_{j=1}^\infty$ be a sequence in $S_r$ that converges to $m$; we may assume without loss of generality that $(m_j)_{j=1}^\infty$ is strictly increasing.  By assumption (A4), there exists an interval $I$ containing $m$ on which $A$ is convex.  This interval contains $[m_{j-1}, m_j]$ for large enough~$j$.  We know that $A(m_{j-1}) = A^{**}(m_{j-1})$, $A(m_j) = A^{**}(m_j)$, and $A$ is convex on $[m_{j-1}, m_j]$; therefore $A = A^{**}$ on $[m_{j-1}, m_j]$, whence $m_j\in \Sigma_0$.  This is impossible, since $\Sigma_0\cap S_r = \emptyset$.  
	\end{proof}
	
	This lemma has the following obvious but useful consequence:
	\begin{COROL}
		The set $\Sigma_+\cup \Sigma_0$ is a union of half-open intervals of the form $(m',m'']$.  In particular, if $m\in \Sigma_+\cup \Sigma_0$, then there exists $m'<m$ such that $A = A^{**}$ on $[m',m]$. 
	\end{COROL}

	\subsection{The sticky particle Cucker--Smale dynamics}
	The entropic solution to the Euler-alignment system \eqref{e:EA} is compatible with the Cucker--Smale system of ODE's \cite{cucker2007emergent}
	\begin{equation}
	\label{eq:CS}
	\frac{\dd x_i}{\dd t} = v_i,\quad
	\frac{\dd v_i}{\dd t} =  \sum_{j=1}^N m_j \phi(x_j - x_i) (v_j - v_i),
	\qquad 
	i = 1, \ldots, N.
	\end{equation}
	subject to a \textit{sticky particle collision rule} (described below).  The masses $m_i$ are fixed; the positions $x_i$ and velocities $v_i$ satisfy some initial conditions
	\[
	x_i(0) = x_i^0, \quad v_i(0) = v_i^0,
	\qquad 
	i = 1, \ldots, N,
	\]
	with
	\[
	x_1^0\le x_2^0\le \cdots \le x_N^0.
	\]
	We do not require the $x_i^0$'s to be distinct, but we insist that \begin{equation}
	\label{e:samevel}
	v_i^0 = v_j^0 \text{ whenever } x_i^0 = x_j^0.  
	\end{equation}
	
	We now specify the collision rules. Define the index cluster $J_i(t)$ to be the collection of indices associated to the agents which are stuck to agent $i$ at time $t$; let $i_*(t)$ and $i^*(t)$ denote the minimum and maximum of $J_i(t)$, respectively: 
	\[
	J_i(t) = \{j:x_j(t) = x_i(t)\} = \{i_*(t), i_*(t)+1, \ldots, i^*(t)\}.
	\]
	A \textit{collision} occurs when $J_i(t)$ changes cardinality. We impose `sticky particle' collision rules as follows:
	\begin{itemize}
		\item Each collision is \emph{completely inelastic}, and agents stick to each other after collisions:
		\begin{equation}
		\label{e:stick}
		J_i(t) \supseteq  J_i(s), \quad \text{ whenever }\,\, t\ge s\ge 0;
		\end{equation}
		\item Collisions conserve momentum:
		\begin{equation}
		\label{e:momentumconscollision}
		v_i(t) = \frac{\sum_{j\in J_i(t)} m_j v_j(t-)}{\sum_{j\in J_i(t)} m_j}.
		\end{equation}
	\end{itemize}
	For convenience, we also assume that the velocities are right continuous, i.e., $v_i(t) = v_i(t+)$.  We also frequently use standard notation for time derivatives to indicate derivatives \textit{from the right}, which always exist under our conventions, even though $v_i(t)$ may experience jump discontinuities.
	
	Next, we introduce an important quantity 
	\begin{equation}\label{eq:psi}
	\psi_i(t) = v_i(t) + \sum_{j=1}^N m_j \Phi(x_i(t) - x_j(t)),
	\qquad 
	i = 1, \ldots, N.
	\end{equation}
	It is conserved in time, in the sense described by the following Lemma.
	\begin{LEMMA}[Conservation of $\psi_i$]\label{lem:psiconservation}
		For any non-collision time $t$, we have 
		\begin{equation}
		\label{e:psiidot=0}
		\frac{\dd}{\dt}\psi_i(t) = 0,\qquad\forall~i=\{1,\ldots,N\}.
		\end{equation}
		For any collision time $t$ (and in fact, any time), we have
		\begin{equation}
		\label{e:psicollision}
		\psi_i(t) = \frac{\sum_{j\in J_i(t)} m_j \psi_j(t-)}{\sum_{j\in J_i(t)} m_j}
		= \frac{\sum_{j\in J_i(t)} m_j \psi_j^0}{\sum_{j\in J_i(t)} m_j},
		\qquad\forall~i=\{1,\ldots,N\}.
		\end{equation}
	\end{LEMMA}
	
	The following lemma describes the behavior at a collision. 
	\begin{LEMMA}[Barycentric lemma]\label{l:barycentric}
		For any $i\in \{1, \ldots, N\}$ and any $t\ge 0$, we have
		\begin{equation}\label{e:barycentric}
		\frac{\sum_{\ell=i_*(t)}^j m_\ell \psi_\ell(t-)}{\sum_{\ell = i_*(t)}^j m_\ell}
		\ge \psi_i(t) \ge 
		\frac{\sum_{\ell=j+1}^{i^*(t)} m_\ell \psi_\ell(t-)}{\sum_{\ell = j+1}^{i^*(t)} m_\ell}, 
		\quad \forall~j= \{i_*(t),\cdots,i^*(t)\}.
		\end{equation}
	\end{LEMMA}
	
	A different version of the barycentric Lemma was previously used by Brenier and Grenier \cite{brenier1998sticky,grenier1995existence} to analyze the 1D pressureless Euler equations.  The latter can be recovered as a special case of \eqref{e:EA} when $\phi\equiv0$, in which case we also have $\psi_i=v_i$. The inequality \eqref{e:barycentric} then means that when collision occurs, the average velocity of the left group of particles has to be larger than the average velocity of the right group of particles. As noticed by the authors in \cite{leslie2023sticky}, the barycentric lemma extends to the case of general (locally integrable) communication protocols $\phi$; this extension is the statement recorded in Lemma \ref{l:barycentric}.
	
	The properties detailed in Lemmas \ref{lem:psiconservation} and \ref{l:barycentric} endow the quantities $\psi_i$ with crucial information about collisions and cluster formation for the discrete sticky particle Cucker--Smale system. The consequences of these properties will be thoroughly investigated in this paper.  
	
	Before moving on, we pause to record a simple identity we will use repeatedly:
	\begin{equation}
	\label{e:xdiff}
	\begin{split}
	\frac{\dd}{\dt} (x_j(t) - x_i(t)) 
	& = \psi_j(t) - \psi_i(t) - \sum_{\ell=1}^N m_\ell \int_{x_i(t)}^{x_j(t)} \phi(y - x_\ell(t))\dy,
	\qquad \forall i,j\in \{1, \ldots, N\}.
	\end{split} 
	\end{equation}

	\subsection{Atomic solutions of the Euler-alignment system}
	\label{sec:atomic}
	We recall the following connection between the Euler-alignment system \eqref{e:EA} and the sticky particle Cucker--Smale dynamics \eqref{eq:CS}.
	\begin{PROP}[{\cite{leslie2023sticky}}]
		Consider the Euler-alignment system \eqref{e:EA} with atomic initial data
		\begin{equation}\label{eq:rhoPN0}
		\rho_N^0(x) = \sum_{i=1}^{N} m_{i,N} \d(x-x_{i,N}^0),
		\qquad 
		P_N^0(x) := \rho_N^0 u_N^0(x) = \sum_{i=1}^{N} m_{i,N} v_{i,N}^0 \d(x-x_{i,N}^0).
		\end{equation}
		There exists a unique entropic solution 
		\begin{equation}\label{eq:rhoPNt}
		\rho_N(x,t) = \sum_{i=1}^{N} m_{i,N} \d(x-x_{i,N}(t)),
		\quad 
		P_N(x,t) = \rho_N u_N(x,t) = \sum_{i=1}^{N} m_{i,N} v_{i,N}(t) \d(x-x_{i,N}(t)),
		\end{equation}
		where $(x_{i,N}(t),v_{i,N}(t))_{i=1}^N$ is the solution to the sticky particle Cucker--Smale dynamics \eqref{eq:CS} with initial data $(m_{i,N},x_{i,N}^0,v_{i,N}^0)_{i=1}^N$.
	\end{PROP}
	For initial data of the form \eqref{eq:rhoPN0}, the cumulative distribution function $M_N^0$ and its generalized inverse $X_N^0$ (defined as in \eqref{eq:M0}) are piecewise constant functions. We write out their formulas presently.  Define  
	\begin{equation}
	\label{e:thiN}
	\th_{i,N} = -\frac12+\sum_{j=1}^i m_{j,N},\quad i=0,\ldots,N,
	\end{equation}
	so that
	\begin{equation}
	m_{i,N} = \th_{i,N} - \th_{i-1,N},\quad i=1,\ldots,N.
	\end{equation}
	Then the corresponding initial data for the scalar balance law \eqref{e:SB}, and its generalized inverse, are
	\begin{equation}
	M_N^0(x) = -\frac12 + \sum_{i=1}^N m_{i,N} \one_{[0,\infty)}(x-x_{i,N}^0) = -\frac12 \one_{(-\infty,x_{1,N}^0)}(x) + \sum_{i=1}^N \th_{i,N} \one_{[x^0_{i,N}, x_{i+1,N}^0)}(x),
	\end{equation}
	\begin{equation}
	X_N^0(m) = \sum_{i=1}^N x_{i,N}^0 \one_{(\th_{i-1,N}, \th_{i,N}]}(m).
	\end{equation}
	Here, we use the convention 
	$x_{N+1,N}^0=+\infty$.  From Definition \ref{def:cluster}, the initial data includes a cluster of size $\sum_{j\in J_i(0)} m_{j,N}$ located at each $x_{i,N}^0$.  The blue step function in the left subplot of Figure \ref{fig:discrete} shows a typical scenario.
	
	Using the atomic initial data \eqref{eq:rhoPN0} in the formula \eqref{e:defA} for the flux yields a piecewise linear function $A_N$ with the formula
	\begin{equation}
	\label{eq:defAN}
	A_N(m) = \sum_{j=1}^{i-1} m_{j,N} \psi_{j,N}^0 + (m-\th_{i-1,N}) \psi_{i,N}^0, 
	\qquad \th_{i-1,N}<m\le \th_{i,N},
	\qquad i=1, \ldots, N,
	\end{equation}	
	for $m\in (-\frac12, \frac12]$, and $A_N(-\frac12) = 0$.  Here the $\psi_{i,N}^0$'s are defined by  
	\begin{equation}
	\psi_{i,N}^0 
	= v_{i,N}^0 + \sum_{j=1}^N m_{j,N} \Phi(x_{i,N}^0 - x_{j,N}^0),
	\qquad i=1, \ldots, N.
	\end{equation}
	We also define $\psi_{i,N}(t)$ analogously to \eqref{eq:psi}:
	\begin{equation}
	\label{e:psiiN}
	\psi_{i,N}(t) := v_{i,N}(t) + \sum_{j=1}^N m_{j,N} \Phi(x_{i,N}(t) - x_{j,N}(t)),
	\qquad 
	i = 1, \ldots, N.
	\end{equation}
	Note that we have 
	\begin{equation}\label{eq:thetaslope}
	m_{i,N} \psi_{i,N}^0 = A_N(\th_{i,N}) - A_N(\th_{i-1,N}),\quad\text{i.e.,}\quad
	\psi_{i,N}^0=\frac{A_N(\th_{i,N}) - A_N(\th_{i-1,N})}{\th_{i,N}-\th_{i-1,N}}.
	\end{equation}
	Hence, $\psi_{i,N}^0$ is the slope of the line segment comprising the graph of $A_N$ in $[\th_{i-1,N},\th_{i,N}]$. For discrete initial data, the flux $A_N$ (together with the masses $(m_{i,N})_{i=1}^N$) therefore carries the information of $(\psi_{i,N}^0)_{i=1}^N$. See the blue piecewise linear curve in the right subplot of Figure~\ref{fig:discrete} for an illustration.  We also note that this process yields initial velocities which satisfy \eqref{e:samevel}, as can be easily checked. 
	
	When collisions occur, we apply \eqref{e:psicollision} and \eqref{eq:thetaslope} to get
	\begin{equation}\label{eq:psiA}
	\psi_{i,N}(t)
	=\frac{\sum_{j\in J_i(t)} m_{j,N} \psi_{j,N}^0}{\sum_{j\in J_i(t)} m_{j,N}}
	=\frac{A_N(\th_{i^*(t),N})-A_N(\th_{i_*(t)-1,N})}{\th_{i^*(t),N}-\th_{i_*(t)-1,N}}.	
	\end{equation}
	In other words, $\psi_{i,N}(t)$ is the slope of secant line through the graph of $A_N$ between the two points $(\th_{i_*(t)-1,N}, A_N(\th_{i_*(t)-1,N}))$ and $(\th_{i^*(t),N}, A_N(\th_{i^*(t),N}))$. For instance, as illustrated in Figure~\ref{fig:discrete}, if agents 2 and 3 stick at time $t$ (and are not stuck to any other agents at time $t$), then 
	\[\psi_{2,N}(t)=\psi_{3,N}(t)=\frac{m_{2,N}\psi_{2,N}^0 + m_{3,N}\psi_{3,N}^0}{m_{2,N} + m_{3,N}}\]
	is the slope of the red secant line through $(\th_{1,N}, A_N(\th_{1,N}))$ and $(\th_{3,N}, A_N(\th_{3,N}))$.
	
	\subsection{The sticky particle approximation}
	\label{sec:approx}
	
	We now reinstate assumptions (A1)--(A4) and consider a \textit{sequence} of atomic solutions $(\rho_N, u_N)_{N=1}^\infty$ that approximate the solution $(\rho, u)$ of interest.  As outlined in the previous subsection, generating an atomic solution $(\rho_N, u_N)$ for any given $N$ amounts to choosing initial data $(m_{i,N}, x_{i,N}^0, v_{i,N}^0)_{i=1}^N$ and running the dynamics for the sticky particle Cucker--Smale system.  We now present our choice of sticky particle initial data and make precise the sense in which the corresponding $(\rho_N, u_N)$ approximates $(\rho, u)$.
	
	We choose our $m_{i,N}$'s in such a way that
	\begin{equation} 
	\label{eq:milimit}
	\tag{D1}
	\lim_{N\to \infty} \max_{1\le i\le N} m_{i,N} = 0.
	\end{equation} 
	Then, we define $\th_{i,N}$ as in \eqref{e:thiN}.  For large enough $N$ (say $N\ge N_0$), we may always choose the $m_{i,N}$'s to satisfy the following additional hypothesis:
	\begin{equation}
	\label{eq:D2}
	\tag{D2}
	S=\p\Sigma_- \subset \{\th_{i,N}\}_{i=0}^N, 
	\qquad N\ge N_0.
	\end{equation}
	We define the remaining sticky particle data in terms of $(M^0, A)$ (which we recall are defined in \eqref{eq:M0}--\eqref{e:defA}), and the $m_{i,N}$'s:
	\begin{equation}
	\label{e:xiN0}
	\tag{D3}
	x_{i,N}^0 = \inf\{x\in \R: M^0(x)\ge \th_{i,N}\} = X^0(\th_{i,N}),
	\end{equation}
	and finally  
	\begin{equation}
	\label{e:psii0Nvi0N}
	\tag{D4}
	\psi_{i,N}^0 = \frac{A(\th_{i,N}) - A(\th_{i-1,N})}{\th_{i,N} - \th_{i-1,N}},
	\qquad v_{i,N}^0 = \psi_{i,N}^0 - \sum_{j=1}^N m_{j,N} \Phi(x_{i,N}^0 - x_{j,N}^0),
	\qquad i=1, \ldots, N.
	\end{equation}
	This gives us all the information we need in order to define $M_N^0$, $X_N^0$, and $A_N$, for each $N\in \N$, using exactly the formulas of the previous subsection.  
	
\begin{figure}[ht]
		\begin{tikzpicture}
		\draw[thick,<->] (6.5,-.3) node[right]{$x$} -- (0,-.3) -- (0,3.7) node[above]{$M^0$};
		\draw[thick,blue] (0,0.2) -- (1,0.2);
		\draw[thick,blue] (1,.8) -- (2.2,.8); 
		\draw[thick,blue] (2.2,1.4) -- (3, 1.4); 
		\draw[thick,blue] (3, 1.8) -- (3.2, 1.8); 
		\draw[thick,blue] (3.2, 2.2) -- (4,2.2); 
		\draw[thick,blue] (4, 2.6) -- (5.5,2.6); 
		\draw[thick,blue] (5.5,3.2) -- (6.5,3.2);
		\node[circle,draw=blue, fill=white, inner sep=0pt,minimum size=3.5pt] at (1,0.2) {};
		\node[circle,draw=blue, fill=blue, inner sep=0pt,minimum size=3.5pt] at (1,0.8) {};
		\node[circle,draw=blue, fill=white, inner sep=0pt,minimum size=3.5pt] at (2.2,.8) {};
		\node[circle,draw=blue, fill=blue, inner sep=0pt,minimum size=3.5pt] at (2.2,1.4) {};
		\node[circle,draw=blue, fill=white, inner sep=0pt,minimum size=3.5pt] at (3, 1.4) {};
		\node[circle,draw=blue, fill=blue, inner sep=0pt,minimum size=3.5pt] at (3, 1.8) {};
		\node[circle,draw=blue, fill=white, inner sep=0pt,minimum size=3.5pt] at (3.2, 1.8) {};
		\node[circle,draw=blue, fill=blue, inner sep=0pt,minimum size=3.5pt] at (3.2, 2.2) {};
		\node[circle,draw=blue, fill=white, inner sep=0pt,minimum size=3.5pt] at (4,2.2) {};
		\node[circle,draw=blue, fill=blue, inner sep=0pt,minimum size=3.5pt] at (4, 2.6) {};
		\node[circle,draw=blue, fill=white, inner sep=0pt,minimum size=3.5pt] at (5.5,2.6) {};
		\node[circle,draw=blue, fill=blue, inner sep=0pt,minimum size=3.5pt] at (5.5,3.2) {};
		
		\draw[dashed] plot[smooth,tension=.6] coordinates{(0.2,0.2) (0.5, 0.3) (1,.8) (2.2, 1.4) (3, 1.8) (3.2, 2.2) (4,2.6) (5.1,3.1) (5.5,3.2)};
		\draw (.2,.2) -- (0,.2) node[left]{$\th_{0}=-\frac12$};
		\draw (.2,.8) -- (0,.8) node[left]{$\th_{1}$};
		\draw (.2,1.4) -- (0,1.4);	
		\draw (.2,1.8) -- (0,1.8) node[left]{$\vdots$\,\,};	
		\draw (.2,2.2) -- (0,2.2);	
		\draw (.2,2.6) -- (0,2.6) node[left]{$\th_{N-1}$};	
		\draw (.2,3.2) -- (0,3.2) node[left]{$\th_{N}=\frac12$};
		
		\draw (1,-.1) -- (1,-.3) node[below]{$x^0_{1}$};
		\draw (2.2,-.1) -- (2.2,-.3) node[below]{$x^0_{2}$};
		\draw (3,-.1) -- (3,-.3) node[below, yshift=-4]{$\cdots$};
		\draw (3.2,-.1) -- (3.2,-.3);
		\draw (4,-.1) -- (4,-.3) node[below]{$x^0_{N-1}$};
		\draw (5.5,-.1) -- (5.5,-.3) node[below]{$x^0_{N}$};
		
		\draw[thick,<->] (14.5,-.3) node[right]{$m$} -- (7.5,-.3) -- (7.5,3.7) node[above]{$A$};
		\draw[dashed] plot[smooth,tension=1] coordinates{(8,3) (9.2,.8) (10.4,1.6) (11.2,1) (12,2.6) (12.8,1.8) (14,3.6)};
		\draw[blue,thick] plot[mark=*] coordinates{(8,3) (9.4,.73) (10.4,1.6) (11.1,.96) (12,2.6) (12.8,1.8) (14,3.6)};
		\draw[red] (9.4,.71) -- (11.1,0.96);

		\draw (8,-.1) -- (8,-.3) node[below]{$\th_{0}$};
		\draw (9.4,-.1) -- (9.4,-.3) node[below]{$\th_{1}$};
		\draw (10.4,-.1) -- (10.4,-.3) node[below]{$\th_{2}$};
		\draw (11.1,-.1) -- (11.1,-.3) node[below]{$\th_{3}$};
		\draw (12,-.1) -- (12,-.3)node[below, yshift=-3]{$\cdots$};
		\draw (12.8,-.1) -- (12.8,-.3) node[below]{$\th_{N-1}$};
		\draw (14,-.1) -- (14,-.3) node[below]{$\th_{N}$};
		
		\draw (7.7,3) -- (7.5,3) node[left]{$0$};
		
		\draw[blue] node at (8.8,2.4) {$\psi_{1}^0$};
		\draw[blue] node at (9.6,1.5) {$\psi_{2}^0$};
		\draw[blue] node at (10.9,1.8) {$\psi_{3}^0$};
		\draw[blue] node at (13.1,3) {$\psi_{N}^0$};
		\draw[red] node at (10.2,.5) {$\psi_{2}(t)$};
		\end{tikzpicture}
		\caption{An illustration of the discretization and the flux, with $N=6$.  We have dropped a subscript $N$ on most of the discretized quantities.  Left: $M^0$ (dashed) and its piecewise constant discretization $M_N$ (solid). Right: The flux $A$ (dashed) and its piecewise linear discretization $A_N$ (solid).  The slope of $A$ between the breakpoints $\th_{j-1}$ and $\th_j$ determines the value of $\psi_j$.  If agents $j$ and $j+1$ collide, the values of $\psi_j(t)$ and $\psi_{j+1}(t)$ are adjusted accordingly; for example, the figure shows the value of $\psi_2(t) = \psi_3(t)$ after agents 2 and 3 have collided.}
		\label{fig:discrete}
	\end{figure}
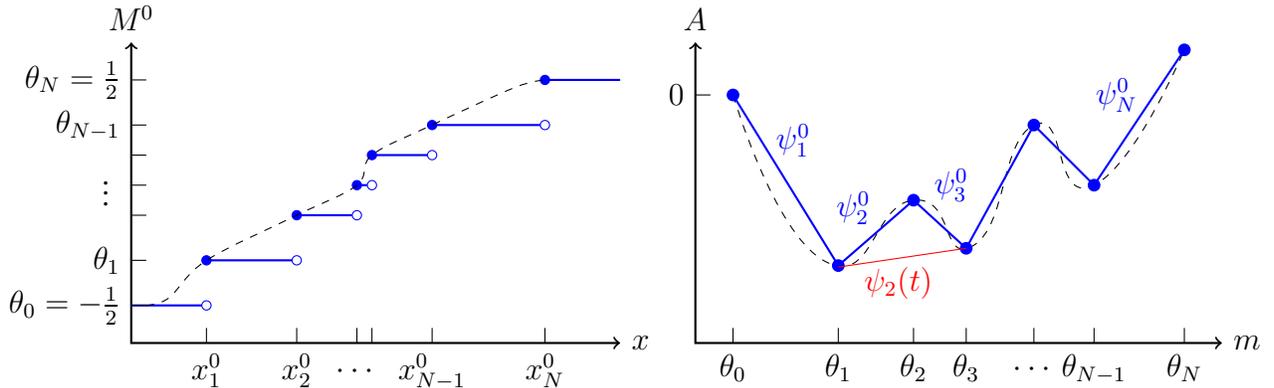

	
	Our discretization scheme has several crucial properties.  First of all, \eqref{eq:milimit} and \eqref{e:xiN0} guarantee that 
	\begin{equation}\label{e:MN0L1conv}
	\|M_N^0 - M^0\|_{L^1(\R)} = \|X_N^0 - X^0\|_{L^1(-\frac12, \frac12)}\to 0, 
	\qquad \text{ as } N\to \infty. 
	\end{equation}
	Next, \eqref{e:psii0Nvi0N} guarantees that $A_N$ (defined in \eqref{eq:defAN}) agrees with $A$ at the breakpoints $\th_{i,N}$ and is linear in between:
	\begin{equation} 
	\label{e:ANvA}
	A_N(\th_{i,N}) =A(\th_{i,N}), \; i=0, \ldots, N; 
	\qquad 
	\text{ and } A_N \text{ is linear on } [\th_{i-1}, \th_i],\;\;  \forall i \in \{1, \ldots, N\}.
	\end{equation} 
	Therefore, $A_N\to A$ uniformly as $N\to \infty$.  See Figure \ref{fig:discrete} for an illustration of the approximation: For initial data $(M^0,A)$ indicated by the dashed curves, the solid blue graphs of $(M_N^0, A_N)$ serve as approximations.

	As the authors proved in \cite{leslie2023sticky}, the entropy solutions generated by discretized initial data satisfying \eqref{eq:milimit}, \eqref{e:xiN0}, \eqref{e:psii0Nvi0N} are a good approximation of the true solution.  Letting $(x_{i,N}(t), v_{i,N}(t))_{i=1}^N$ denote the sticky particle Cucker--Smale dynamics associated to the initial data $(m_{i,N}, x_{i,N}^0, v_{i,N}^0)_{i=1}^N$, the entropy solution of \eqref{e:SB} associated to the initial data $M_N^0$ and flux $A_N$ defined above is 
	\[
	M_N(x,t) = -\frac12 + \sum_{i=1}^N m_{i,N} \one_{[0,\infty)}(x-x_{i,N}(t)).
	\]
	Its generalized inverse is
	\[
	X_N(m,t) = \sum_{i=1}^N x_{i,N}(t) \one_{(\th_{i-1,N}, \th_{i,N}]}(m).
	\]
	We have the following statement from \cite{leslie2023sticky}:
	\begin{PROP}\label{prop:stability}
		For any fixed $t\geq0$, the solutions $M_N(\cdot,t)$ approximate $M(\cdot,t)$ in the following sense:
		\begin{equation}
		\label{e:MNL1conv}
		M_N(\cdot,t) - M(\cdot,t) \to 0 \,\, \text{ in } L^1(\R). 
		\end{equation}	
	\end{PROP}
	The convergence \eqref{e:MNL1conv} is equivalent to saying that $(\rho_N(t))_{N=1}^\infty$ converges to $\rho(t)$ in the Wasserstein-1 metric, and it was also proved in \cite{leslie2023sticky} that $(\rho_N u_N(t))_{N=1}^\infty$ converges weak-$*$ in the sense of measures to $\rho u(t)$.  The convergence \eqref{e:MNL1conv} is what we use below.  More precisely, we note that since
	\[
	\|X_N(\cdot,t) - X(\cdot,t)\|_{L^1(-\frac12, \frac12)} = \|M_N(\cdot,t) - M(\cdot,t)\|_{L^1(\R)},
	\]
	the convergence \eqref{e:MNL1conv} implies the existence of a subsequence $(X_{N_k}(\cdot,t))_{k=1}^\infty$ such that
	\begin{equation}
	\label{e:XNkconv}
	X_{N_k}(m,t)\to X(m,t) \text{ as } k\to \infty,\quad \text{for almost every } m\in(-\tfrac12, \tfrac12].
	\end{equation}
	We will be able to leverage this almost-everywhere convergence, combined with the monotonicity of $X(\cdot,t)$, to analyze the formation of clusters in the sense of Definition \ref{def:cluster}.
	
	We end this subsection by discussing the significance of the requirement \eqref{eq:D2} in our approximation scheme.  It guarantees that 
	\begin{equation}
	\label{e:supersame}
	\{m:A(m)>A^{**}(m)\} = \{m:A_N(m)>A_N^{**}(m)\},
	\end{equation}
	so that the supercritical region of mass labels associated to $(\rho_N, u_N)$ is the same as the one for $(\rho, u)$.  
	
	On account of \eqref{e:ANvA} and \eqref{e:supersame}, it becomes superfluous for our purposes to set additional notation for the discretized flux; we write our analysis below in terms of $A$ only.
	
	\subsection{A maximum principle}
	
	We record one final preliminary statement here; it is a uniform-in-$N$ maximum principle for the velocities associated to our discretization.  
	
	\begin{LEMMA}
		\label{l:maxprinc}
		Given a sequence of discretizations satisfying \eqref{eq:milimit}--\eqref{e:psii0Nvi0N}, there exists an $N$-independent positive constant $u_{\max}$ such that
		\begin{equation}
		\label{e:maxprinc}
		\max_{1\le i\le N} \sup_{t\ge 0} |v_{i,N}(t)| \le \max_{1\le i\le N} |v_{i,N}^0|\le u_{\max},
		\qquad \forall N\in \N.  
		\end{equation}
	\end{LEMMA}	
	We refer to \cite{leslie2023sticky} for the details of the proof.
	
	\section{The subcritical regime} 
	
	\label{sec:sub} 
	
	In this section, we prove part I of Theorem \ref{t:regions} by establishing a more general statement, namely that mass labels from distinct $L(m)$'s can never belong to the same finite- or infinite-time cluster.  Now, it may of course happen that $L(m')\ne L(m'')$ without either of $m'$ or $m''$ belonging to $\Sigma_+$, but such labels are `related subcritically,' so to speak, in that their images remain separated under $X(\cdot, t)$ for essentially the same reasons as do different mass labels in $\Sigma_+$.  As already noted above, $L(m)$ is a singleton for each $m\in \Sigma_+$, so establishing the statement mentioned above proves in particular that there can never be \textit{any} finite- or infinite-time clusters in $\Sigma_+$.  
	
	We begin by proving a discrete version of this more general statement.  The finite-time part (Proposition \ref{p:Lm} below) can be deduced from the barycentric lemma only---no explicit reference to the equations governing the dynamics is required.  We bootstrap our finite-time statement into an infinite-time one with the aid of \eqref{e:xdiff}, proving that the distance between agents associated to different $L(m)$'s does not tend to zero, and thus that agents from different $L(m)$'s cannot be part of the same \textit{infinite}-time cluster either.  Lemma \ref{l:subdisc*} gives a time-independent lower bound on the distance between such agents; we use it in the proof of the continuum version of our statement, which we establish in Section \ref{ssec:subcont}.
	
	In this section, we work with a fixed entropically selected solution $(\rho, u)$ of \eqref{e:EA} associated to initial data $(\rho^0, u^0)$, and we assume throughout that (A1)--(A4) hold.  
	
	\subsection{The discrete setting}
	
	In this subsection, we fix $N\in \N$ large and a discretization $(m_{i,N},x_{i,N}^0,v_{i,N}^0)_{i=1}^N$ of $(\rho^0, u^0)$ satisfying \eqref{eq:D2}--\eqref{e:psii0Nvi0N}. We let $(x_{i,N}(t), v_{i,N}(t))_{i=1}^N$ denote the associated sticky particle Cucker--Smale dynamics.  
	
	Our first main statement of this subsection is the following proposition.
	\begin{PROP}
		\label{p:Lm}
		Fix $i\in \{1, \ldots, N\}$ and $t\ge 0$. If $i_*(t)<i^*(t)$, then 
		$A^{**} \text{ is linear on } [\th_{i_*(t)-1,N}, \th_{i^*(t),N}]$.
	\end{PROP}
	This proposition says that, at the discrete level, finite-time collisions between agents are confined to a single $L(m)$.  In particular, no discrete finite-time clustering can occur anywhere in $\Sigma_+$.  Most of the work involved in proving the proposition goes into establishing the case considered in the following lemma.  The slightly stronger conclusion available in this case will also be useful later in Section \ref{sec:super}.
	
	\begin{LEMMA}
		\label{l:inductionlemma}
		If $i<i^*(t)$ and $\th_{i,N}\notin \Sigma_-$, then $\th_{i_*(t)-1,N}, \th_{i^*(t),N}\notin \Sigma_-$ as well, i.e.,
		\begin{equation}
		\label{e:inductionlemma}
		A(\th_{i_*(t)-1,N}) = A^{**}(\th_{i_*(t)-1,N}) 
		\qquad \text{ and }  \qquad A(\th_{i^*(t),N}) = A^{**}(\th_{i^*(t),N}).
		\end{equation}
		Consequently, we have that
		\begin{equation}
		\label{e:A**lindisc}
		A^{**} \text{ is linear on } [\th_{i_*(t)-1,N}, \th_{i^*(t),N}].
		\end{equation}
	\end{LEMMA}
	The `consequently' claim here follows from Lemma \ref{l:A**segments}.  Note that $i<i^*(t)$ implies that $X_N(\cdot, t)$ is constant on some interval containing $(\th_{i-1,N}, \th_{i+1,N}]$.  In particular, $\th_{i,N}$ lies in the interior of this interval.  This is why we assume $i<i^*(t)$ rather than $i_*(t)<i^*(t)$ here.
	
	Using Lemma \ref{l:inductionlemma}, we give the short proof of the rest of Proposition \ref{p:Lm}.  Then we prove the lemma.  
	
	\begin{proof}[Proof of Proposition \ref{p:Lm}]
		
		If $\th_{j,N}\in \Sigma_-$ for all $j\in \{i_*(t), \ldots, i^*(t)-1\}$, then assumption \eqref{eq:D2} implies that $(\th_{i_*(t)-1,N}, \th_{i^*(t),N})$ must lie in a single connected component of $\Sigma_-$, on which $A^{**}$ is linear, so we are done.  Otherwise, we have $\th_{j,N}\notin \Sigma_-$ for some $j\in \{i_*(t), \ldots, i^*(t)-1\}$.  Then since $j<i^*(t) = j^*(t)$,  the conclusion of Lemma \ref{l:inductionlemma} holds, with $i$ replaced by $j$.  But since agents $i$ and $j$ have collided at time~$t$, the intervals $[\th_{i_*(t)-1,N}, \th_{i^*(t),N}]$ and $[\th_{j_*(t)-1,N}, \th_{j^*(t),N}]$ are the same.  This completes the proof.
	\end{proof}

	\begin{proof}[Proof of Lemma \ref{l:inductionlemma}]
		
		Let $t_1$ denote the first time  satisfying $i < i^*(t_1)$, and let $t_2,\ldots, t_m$ denote any subsequent times (if any) where agent $i$ is involved in a collision.  We prove that \eqref{e:inductionlemma} holds at each time~$t_n$, $n=1, \ldots, m$, which will prove that \eqref{e:inductionlemma} is valid on $[t_1, t_2), [t_2, t_3),\ldots, [t_m,+\infty)$ and thus at any time $t$ such that $i<i^*(t)$.  We argue inductively on $n\in \{1, \ldots, m\}$, splitting our `base case' into two subcases, namely $t_1 = 0$ and $t_1>0$. 
		
		If $t_1 = 0$, then our discretization procedure ensures that $\psi_{i_*(0),N}^0 =\cdots = \psi_{i^*(0),N}^0$, which guarantees that the points $(\th_{i_*(0)-1,N}, A(\th_{i_*(0)-1,N})$, $(\th_{i,N}, A(\th_{i,N}))$, and $(\th_{i^*(0),N}, A(\th_{i^*(0),N})$ are collinear.  Remembering that $A(\th_{i,N}) = A^{**}(\th_{i,N})$ and applying Lemma \ref{l:A**segments} implies that \eqref{e:inductionlemma} holds at time $t_1 = 0$.  
		
		If $t_1>0$, then
		\[
		\frac{\sum_{j={i+1}}^{i^*(t_1)} m_{j,N} \psi_{j,N}(t_1-)}{\sum_{j=i+1}^{i^*(t_1)} m_{j,N}} 
		= \frac{\sum_{j={i+1}}^{i^*(t_1)} m_{j,N} \psi_{j,N}^0}{\sum_{j=i+1}^{i^*(t_1)} m_{j,N}} = \frac{A(\th_{i^*(t_1),N}) - A(\th_{i,N})}{\th_{i^*(t_1),N} - \th_{i,N}}
		\]
		and similarly
		\[
		\frac{\sum_{j=i_*(t_1)}^{i} m_{j,N} \psi_{j,N}(t_1-)}{\sum_{j=i_*(t_1)}^{i} m_{j,N}} = \frac{A(\th_{i,N}) - A(\th_{i_*(t_1)-1,N})}{\th_{i,N} - \th_{i_*(t_1)-1,N}}.
		\]
		The barycentric lemma then tells us that 
		\[
		\frac{A(\th_{i^*(t_1),N}) - A(\th_{i,N})}{\th_{i^*(t_1),N} - \th_{i,N}}
		\le 
		\frac{A(\th_{i,N}) - A(\th_{i_*(t_1)-1,N})}{\th_{i,N} - \th_{i_*(t_1)-1,N}}.
		\]
		On the other hand, our assumption that $A(\th_{i,N}) = A^{**}(\th_{i,N})$, the fact that $A\ge A^{**}$ in general, and Lemma \ref{l:convsec} combine to give us 
		\begin{align*}
		\frac{A(\th_{i^*(t_1),N}) - A(\th_{i,N})}{\th_{i^*(t_1),N} - \th_{i,N}}
		& \ge \frac{A^{**}(\th_{i^*(t_1),N}) - A^{**}(\th_{i,N})}{\th_{i^*(t_1),N} - \th_{i,N}} \\
		& \ge 
		\frac{A^{**}(\th_{i,N}) - A^{**}(\th_{i_*(t_1)-1,N})}{\th_{i,N} - \th_{i_*(t_1)-1,N}}
		\ge 
		\frac{A(\th_{i,N}) - A(\th_{i_*(t_1)-1,N})}{\th_{i,N} - \th_{i_*(t_1)-1,N}}.
		\end{align*}
		The inequalities here thus must actually be equalities.  This forces  \eqref{e:inductionlemma} to hold at time $t_1$.  
		
		Assume inductively that \eqref{e:inductionlemma} holds at time $t_n$ for some $n<m$.  Denote $j = i_*(t_n)-1$.  If agent $i$ collides with an agent on its left at time $t_{n+1}$, then $t_{n+1}$ is the first time when agent $j$ collides with agent $j+1 = i_*(t_n)$, and furthermore, we have $A(\th_{j,N}) = A^{**}(\th_{j,N})$ by our inductive assumption.  
		
		Applying the conclusion from our `base case' to agent $j$, we conclude that $A^{**}$ is linear on the interval $[\th_{j_*(t_{n+1})-1,N}, \th_{j^*(t_{n+1}),N}]$, which is the same as $[\th_{i_*(t_{n+1})-1,N}, \th_{i^*(t_{n+1}),N}]$, since agents $i$ and $j$ collide at time $t_{n+1}$.  Furthermore, 
		\[
		A(\th_{i_*(t_{n+1})-1,N}) 
		= A(\th_{j_*(t_{n+1})-1,N}) 
		= A^{**}(\th_{j_*(t_{n+1})-1,N})
		= A^{**}(\th_{i_*(t_{n+1})-1,N}),
		\]
		and similarly $A(\th_{i^*(t_{n+1}),N}) = A^{**}(\th_{i^*(t_{n+1}),N})$.
		We have thus proved the inductive step in the case where agent $i$ experiences a collision from its left at time $t_{n+1}$.  The case of a collision from the right is entirely similar.  \end{proof}
	
	We now turn our attention to (the absence of) infinite-time clustering, for which we must rely on equation \eqref{e:xdiff}.  The following corollary of Lemma \ref{l:inductionlemma} will help us control the difference $\psi_{j,N}(t) - \psi_{i,N}(t)$ when $\th_{i,N}$ and $\th_{j,N}$ belong to different $L(m)$'s (and in particular when they lie in $\Sigma_+$).
	
	\begin{COROL}
		\label{c:Apinned}
		Assume that $\th_{i,N}\notin \Sigma_-$.  Then
		\begin{equation}
		\label{e:Apinnedi}
		\psi_{i,N}(t)\le \frac{A^{**}(\th_{i,N}) - A^{**}(\th_{i-1,N})}{\th_{i,N} - \th_{i-1,N}},
		\qquad \forall t\ge 0,\; 
		\text{ if } i\in \{1, \ldots, N\};
		\end{equation}		
		\begin{equation}
		\label{e:Apinnedi+1}
		\psi_{i+1,N}(t) \ge \frac{A^{**}(\th_{i+1,N}) - A^{**}(\th_{i,N})}{\th_{i+1,N} - \th_{i,N}}
		\qquad \forall t\ge 0, \; 
		\text{ if } i\in \{0, \ldots, N-1\}.
		\end{equation}	
	\end{COROL}
	
	\begin{proof}
		We prove only \eqref{e:Apinnedi}; the proof of \eqref{e:Apinnedi+1} is similar.  Choose $t\ge 0$.  If $i^*(t) = i$, then 
		\[
		\psi_{i,N}(t) = \frac{A(\th_{i,N}) - A(\th_{i_*(t)-1,N})}{\th_{i,N} - \th_{i_*(t)-1,N}} 
		\le \frac{A^{**}(\th_{i,N}) - A^{**}(\th_{i_*(t)-1,N})}{\th_{i,N} - \th_{i_*(t)-1,N}} 
		\le \frac{A^{**}(\th_{i,N}) - A^{**}(\th_{i-1,N})}{\th_{i,N} - \th_{i-1,N}}.
		\]
		If $i^*(t)>i$, then the previous lemma implies that $A^{**}$ is linear on $[\th_{i_*(t)-1,N}, \th_{i^*(t),N}]$, which of course contains $[\th_{i-1,N}, \th_{i,N}]$.  Therefore,
		\[
		\psi_{i,N}(t) 
		= \frac{A(\th_{i^*(t),N}) - A(\th_{i_*(t)-1,N})}{\th_{i^*(t),N} - \th_{i_*(t)-1,N}} 
		\le \frac{A^{**}(\th_{i^*(t),N}) - A^{**}(\th_{i_*(t)-1,N})}{\th_{i^*(t),N} - \th_{i_*(t)-1,N}}
		= \frac{A^{**}(\th_{i,N}) - A^{**}(\th_{i-1,N})}{\th_{i,N} - \th_{i-1,N}}.
		\]
		This completes the proof.  
	\end{proof}
	
	The following lemma is our second main statement of this subsection.  It gives a time-independent lower bound on the distance between agents from different $L(m)$'s and in particular shows that such agents cannot belong to the same infinite-time cluster.  
	\begin{LEMMA}
		\label{l:subdisc*}
		Fix $i,j\in \{1, \ldots, N\}$ with $i\le j$.  Assume that  there exists $\s>0$ such that
		\begin{equation}
		\label{e:i<j+1}
		\frac{A^{**}(\th_{j+1,N}) - A^{**}(\th_{j,N})}{\th_{j+1,N} - \th_{j,N}}
		- \frac{A^{**}(\th_{i,N}) - A^{**}(\th_{i-1,N})}{\th_{i,N} - \th_{i-1,N}} \ge 2\s >0.
		\end{equation}	
		Then agents $i$ and $j+1$ never collide.  In fact, we have the lower bound
		\begin{equation}
		\label{e:discretesubbd}
		x_{j+1,N}(t) - x_{i,N}(t) \ge \max\bigg\{ x_{j+1,N}^0-x_{i,N}^0 - t u_{\max},  \min\{t\s, \eta\}\bigg\} \ge c>0,	
		\qquad \forall t\ge 0,
		\end{equation}
		where $u_{\max}$ is as in Lemma \ref{l:maxprinc} and $\eta>0$ is chosen so that $|\int_z^w \phi(r)\dr|<\s$ whenever $|z - w|<\eta$.  
	\end{LEMMA}
	
	\begin{proof}[Proof of Lemma \ref{l:subdisc*}]
		
		The fact that agents $i$ and $j+1$ do not collide in finite time is a consequence of Proposition \ref{p:Lm}, since \eqref{e:i<j+1} implies that $A^{**}$ is not linear on any interval containing both $\th_{i-1,N}$ and $\th_{j,N}$.  We need only prove the lower bound.  In fact, it is clear that 
		\[
		x_{j+1,N}(t) - x_{i,N}(t)\ge x_{j+1,N}^0 - x_{i,N}^0 - t u_{\max}.
		\]
		We therefore concentrate on proving that 
		\begin{equation}
		\label{e:minbddisc}
		x_{j+1,N}(t) - x_{i,N}(t)\ge \min\{t\s, \eta\}.  
		\end{equation}
		
		We first give the proof under the additional assumption that $\th_{i,N}, \th_{j,N}\notin \Sigma_-$.  In this case, Corollary~\ref{c:Apinned} tells us that 
		\[
		\psi_{j+1,N}(t) - \psi_{i,N}(t) \ge 2\s >0,
		\qquad \forall t\ge 0.
		\] 
		Therefore, if $\t$ is any time such that  $x_{j+1,N}(\t)-x_{i,N}(\t)<\eta$, then the identity \eqref{e:xdiff} and the previous step tell us that
		\begin{align*}
		\frac{\dd}{\ds}(x_{j+1,N}(s)-x_{i,N}(s))\bigg|_{s = \t} & = \psi_{j+1,N}(\t)-\psi_{i,N}(\t)-\sum_{\ell=1}^N m_{\ell,N}\int_{x_{i,N}(\t)}^{x_{j+1,N}(\t)} \phi(y - x_{\ell, N}(\t)) \dy \ge \s.
		\end{align*}
		It follows immediately that 
		\begin{equation}
		\label{e:discretesubbd2}
		x_{j+1,N}(t) - x_{i,N}(t) \ge \min\{ x_{j+1,N}^0-x_{i,N}^0 + t \s, \eta\},
		\qquad \forall t\ge 0,
		\end{equation}	
		which in particular implies \eqref{e:minbddisc}.  
		
		In the more general setting where one or both of $\th_{i,N}, \th_{j,N}$ may belong to $\Sigma_-$, we replace $i$ and $j$ with the closest indices $I$ and $J$ satisfying $i\le I\le J\le j$ and $\th_{I,N}, \th_{J,N}\notin \Sigma_-$.  This is possible since $\th_{i,N}$ and $\th_{j,N}$ cannot belong to the same connected component of $\Sigma_-$; otherwise $A^{**}$ would be linear on $[\th_{i-1,N}, \th_{j,N}]$, contrary to our assumption.  
		
		More specifically, if $\th_{i,N}\in \Sigma_-$, we choose $I$ so that $\th_{I,N}$ is the right endpoint of the connected component of $\Sigma_-$ to which $\th_{i,N}$ belongs.  If $\th_{j,N}\in \Sigma_-$, then we choose $J$ so that $\th_{J,N}$ is the corresponding left endpoint.  The point is that now we have $\th_{I,N},\th_{J,N}\in \Sigma_-$ and 
		\[
		\frac{A^{**}(\th_{J+1,N}) - A^{**}(\th_{J,N})}{\th_{J+1,N} - \th_{J,N}}
		- \frac{A^{**}(\th_{I,N}) - A^{**}(\th_{I-1,N})}{\th_{I,N} - \th_{I-1,N}}
		= 
		\frac{A^{**}(\th_{j+1,N}) - A^{**}(\th_{j,N})}{\th_{j+1,N} - \th_{j,N}}
		- \frac{A^{**}(\th_{i,N}) - A^{**}(\th_{i-1,N})}{\th_{i,N} - \th_{i-1,N}},
		\]
		so that, applying the logic of the previous case, we have 
		\[
		x_{j+1,N}(t) - x_{i,N}(t)\ge x_{J+1,N}(t) - x_{I,N}(t)\ge \min\{x^0_{J+1,N} - x^0_{I,N} + t\s, \eta\}\ge \min\{t\s, \eta\},
		\]
		as needed. This completes the proof.
	\end{proof}
	
	\begin{REMARK}
		If $\th_{i,N}, \th_{j,N}\in \Sigma_+$ and $i\le j$, then the hypotheses of Lemma \ref{l:subdisc*} are satisfied automatically for some $\s>0$. Indeed, if the left side of \eqref{e:i<j+1} were equal to zero, it would force $A^{**}$ to be linear on $[\th_{i-1,N}, \th_{j+1,N}]$ contradicting the definition of $\Sigma_+$.  
	\end{REMARK}

	\subsection{The continuum setting}
	
	\label{ssec:subcont}
	
	We are now ready to prove the full continuum version of the statement that mass labels from distinct $L(m)$'s never belong to the same cluster.  We first show in Lemma \ref{l:mL} that this is true at time zero; then we give the full statement in Theorem \ref{thm:sub} below.  
	
	\begin{LEMMA}
		\label{l:mL}
		Let $(\rho, u)$ and all related notation be defined as in assumptions (A1)--(A4).  Fix $m''\in~(-\frac12, \frac12]$ and define ${m}_L'' = \inf L({m}'')$. Then 
		\begin{itemize}
			\item [(i)] $A({m}_L'') = A^{**}({m}_L'')$. (I.e., $m_L''\notin \Sigma_-$.)
			\item [(ii)] If ${m}'<{m}_L''$, then $X^0({m}')<X^0({m}'')$.
		\end{itemize}
	\end{LEMMA}
	\begin{proof}
		(i) If  $A({m}_L'') > A^{**}({m}_L'')$, then ${m}_L''\in \Sigma_-$.  Let $(m_-, m_+)$ denote the connected component of $\Sigma_-$ containing ${m}_L''$.  Then $A^{**}$ is linear on the nontrivially overlapping intervals $[m_-, m_+]$ and $[{m}_L'', {m}'']$, hence on their union $[m_-, {m}'']$.  Since $m_- < {m}_L''$, this contradicts the definition of ${m}_L''$.
		
		(ii) We again argue by contradiction.  Suppose $X^0({m}') = X^0({m}'')$; then the definition \eqref{e:defA} of $A$ implies that $A$ is linear on $[{m}',{m}'']$.  In particular, $({m}_L'', A({m}_L'')) = ({m}_L'', A^{**}({m}_L''))$ lies on the segment joining $({m}', A({m}'))$ and $({m}'', A({m}''))$; by Lemma \ref{l:A**segments}, it follows that  $A^{**}$ is linear on $[{m}', {m}'']$, contradicting the definition of ${m}_L''$.	
	\end{proof}

	\begin{THEOREM}\label{thm:sub}
		Let $(\rho, u)$ and all related notation be defined as in assumptions (A1)--(A4).  Fix $m''\in~(-\frac12, \frac12]$. If $m'<\inf L(m'')$, then there exists a time-independent constant $c>0$ such that 
		\begin{equation}
		X(m'', t) - X(m',t) \ge c>0,\quad\forall~t\ge 0.
		\end{equation}
	\end{THEOREM}
	
	\begin{proof}
		Define
		\[
		m_L'':=\inf L(m'').
		\]
		There are two cases to consider, namely $m_L''\in \Sigma_+$ and $m_L''\notin \Sigma_+$.	In either case, we choose a discretization satisfying  \eqref{eq:milimit}--\eqref{e:psii0Nvi0N} \textit{and additionally}, $m_L''\in \{\th_{i,N}\}_{i=0}^N$ for sufficiently large $N$. (This will only be actually used in the second case.) For each such $N$, we choose $\ell$ such that $\th_{\ell,N} = m_L''$.  (This $\ell$ will of course depend on $N$, but we suppress this in the notation.) Proposition \ref{p:Lm} guarantees that agents $\ell$ and $\ell+1$ can never collide. 
		
		We split the remainder of the argument into the two cases mentioned above.
		
		\textbf{Case 1: $m_L''\in \Sigma_+$.} In this case, we may assume without loss of generality that $m'' = m_L''\in \Sigma_+$ and that $A = A^{**}$ on $[m',m'']$.  
		
		\textbf{Step 1.1.} Fix $\widetilde{m}'$ and $\overline{m}$ such that 
		\[
		m'<\widetilde{m}'<\overline{m}<m''.
		\]
		Then by Lemma \ref{l:convsec}, together with the fact that $A$ is \textit{not} linear on $[\widetilde{m}', m'']$, we deduce that
		\[
		\frac{A(m'') - A(\overline{m})}{m''-\overline{m}} 
		> \frac{A(\overline{m}) - A(\widetilde{m}')}{\overline{m} - \widetilde{m}'}.
		\]
		Choose $\widetilde{m}''\in (\overline{m}, m'')$ close enough to $m''$ so that 
		\[
		\frac{A(\widetilde{m}'') - A(\overline{m})}{\widetilde{m}''-\overline{m}} 
		> \frac{A(\overline{m}) - A(\widetilde{m}')}{\overline{m} - \widetilde{m}'},
		\]
		and let $2\s$ denote the difference between the left and right sides:
		\begin{equation}
		\label{e:def2sigma}
		2\s := \frac{A(\widetilde{m}'') - A(\overline{m})}{\widetilde{m}''-\overline{m}}  - \frac{A(\overline{m}) - A(\widetilde{m}')}{\overline{m} - \widetilde{m}'} >0.
		\end{equation} 
		
		\textbf{Step 1.2.} For each $N\in \N$, let $(m_{i,N}, x_{i,N}^0, v_{i,N}^0)_{i=1}^N$ be a discretization of $(\rho^0, u^0)$, and assume that these discretizations satisfy \eqref{eq:milimit}--\eqref{e:psii0Nvi0N}.  Using the notation of Section \ref{sec:preliminaries}, let us fix a time $t>0$ and choose a corresponding subsequence $(N_k)_{k=1}^\infty$ such that the a.e. $\!\!$ convergence \eqref{e:XNkconv} holds at time zero and time $t$. Then we can find $\th'\in (m', \widetilde{m}')$ and $\th''\in (\widetilde{m}'', m'')$, such that 
		\begin{equation}\label{eq:XNnconv}
		\begin{split} 
		X_{N_k}(\th',t)\to X(\th',t),\quad 
		& X_{N_k}(\th'',t)\to X(\th'',t), \\
		X_{N_k}^0(\th')\to X^0(\th'),\quad 
		& X_{N_k}^0(\th'')\to X^0(\th''),
		\end{split} 	
		\qquad \text{ as } k\to \infty.	
		\end{equation} 
		
		For each large enough $N$, we choose indices $i,j$ so that 
		\[
		\th' \in (\th_{i-1,N},\th_{i,N}]\subset (m', \widetilde{m}'),
		\qquad 
		\th'' \in (\th_{j,N},\th_{j+1,N}] \subset (\widetilde{m}'', m'').
		\]
		Note that $i$ and $j$ depend on $N$, but we suppress this dependence in the notation.  The point is that for all sufficiently large $N\in \N$, we have
		\begin{equation}
		\label{e:labeling}
		X_N(\th',s) = x_{i,N}(s), 
		\qquad X_N(\th'',s) = x_{j+1,N}(s),
		\qquad \forall s\ge 0,
		\end{equation}
		so we can track the approximate positions of mass labels $\th'$ and $\th''$ at time $t$ by looking at $x_{i,N}(t)$ and $x_{j+1,N}(t)$, which are simpler to analyze.
		
		\textbf{Step 1.3.}
		Since $A = A^{**}$ is convex on $[m',m'']$, it follows from Lemma \ref{l:convsec} and our choice of $i,j$ that
		\begin{align*}
		\frac{A^{**}(\th_{i,N}) - A^{**}(\th_{i-1},N)}{\th_{i,N} - \th_{i-1,N}} 
		& \le \frac{A(\overline{m})-A(\widetilde{m}')}{\overline{m}-\widetilde{m}'}, & \text{ since } \th_{i-1,N}<\th_{i,N}< \widetilde{m}'\le  \overline{m}, \\
		\frac{A^{**}(\th_{j+1,N}) - A^{**}(\th_{j},N)}{\th_{j+1,N} - \th_{j,N}} 
		& \ge \frac{A(\widetilde{m}'') - A(\overline{m})}{\widetilde{m}'' - \overline{m}}, & \text{ since } \overline{m} < \widetilde{m}'' < \th_{j,N} < \th_{j+1,N}.
		\end{align*}
		Taking the difference of these inequalities and recalling our definition of $\s$ from \eqref{e:def2sigma}, we see that the hypotheses of Lemma \ref{l:subdisc*} are satisfied.  Then, in accordance with Lemma \ref{l:subdisc*}, we have the following lower bound at time $t$:
		\[
		x_{j+1,N}(t) - x_{i,N}(t) \ge \max\bigg\{ x_{j+1,N}^0-x_{i,N}^0 - t u_{\max},  \min\{t\s, \eta\}\bigg\}.
		\]
		(Actually, a slightly better lower bound holds; c.f. the proof of Lemma \ref{l:subdisc*}, but we use the more general version above so that this step generalizes to Case 2.) Substituting  \eqref{e:labeling} into this lower bound yields
		\[
		X_N(\th'', t) - X_N(\th',t)\ge \max\bigg\{ X_N^0(\th'') - X_N^0(\th') - t u_{\max},  \min\{t\s, \eta\}\bigg\}.
		\]
		
		Then taking $N\to \infty$ along the subsequence $(N_k)_{k=1}^\infty$, we conclude that 
		\[
		X(\th'', t) - X(\th', t) \ge \max\bigg\{ X^0(\th'') - X^0(\th') - t u_{\max},  \min\{t\s, \eta\}\bigg\}.
		\]
		
		\textbf{Step 1.4.} Recalling that 
		\[
		m'<\th'<\widetilde{m}'<\widetilde{m}''<\th''<m'',
		\]
		we may conclude from the previous step that 
		\begin{equation}
		\label{e:1.4}
		X(m'', t) - X(m', t) \ge \max\bigg\{ X^0(\widetilde{m}'') - X^0(\widetilde{m}') - t u_{\max},  \min\{t\s, \eta\}\bigg\}.
		\end{equation}
		Since $A^{**}$ is not linear on $[\widetilde{m}', \widetilde{m}'']$ (by \eqref{e:def2sigma}), we have $\widetilde{m}'<\inf L(\widetilde{m}'')$. Lemma \ref{l:mL} then guarantees that $X^0(\widetilde{m}')<X^0(\widetilde{m}'')$, so that the quantity on the right side of \eqref{e:1.4} is strictly positive.  The proof of Case 1 is now complete.  
		
		\textbf{Case 2: $m_L''\notin \Sigma_+$.} In this case, we may assume without loss of generality that $A^{**}$ is linear on $[m',m_L'']$. We note also that we must have $m_L''<m''$, with $A^{**}$ linear on $[m_L'',m'']$.  (Indeed, the only way for $m_L''$ and $m''$ to be equal is if $L(m'') = \{m''\}$, in which case $m_L'' = m''\in \Sigma_+$, contrary to the case we are considering.) Of course, the slopes of $A^{**}$ on $[m',m_L'']$ and $[m_L'',m'']$ are different.  
		
		\textbf{Step 2.1.} As in Case 1, choose $t>0$ and a corresponding subsequence $(N_k)_{k=1}^\infty$; choose $\widetilde{m}'\in (m',m_L'')$, $\widetilde{m}''\in (m_L'',m'')$, and then $\th\in (m', \widetilde{m}')$, $\th''\in (\widetilde{m}'',m'')$ such that \eqref{eq:XNnconv} holds.
		
		For each large enough $N$, we choose indices $i,j$ so that 
		\[
		\th' \in (\th_{i-1,N},\th_{i,N}]\subset (m', \widetilde{m}'),
		\qquad 
		\th' \in (\th_{j,N},\th_{j+1,N}] \subset (\widetilde{m}'', m'').
		\]

		\textbf{Step 2.2.} Recall that we have chosen $\ell\in \N$ to satisfy $\th_{\ell,N} = m_L''$.  Whenever $N$ is chosen large enough so that $m'\le\th_{\ell-1,N}$ and $\th_{\ell+1,N}\le m''$, we must have 
		\[
		\begin{split} 
		& \hspace{-15 mm} \frac{A^{**}(\th_{j+1,N}) -  A^{**}(\th_{j,N})}{\th_{j+1,N} - \th_{j,N}} 
		- \frac{A^{**}(\th_{i,N}) - A^{**}(\th_{i-1,N})}{\th_{i,N} - \th_{i-1,N}} \\
		& \ge 
		\frac{A^{**}(\th_{\ell+1,N}) -  A^{**}(\th_{\ell,N})}{\th_{\ell+1,N} - \th_{\ell,N}} 
		- \frac{A^{**}(\th_{\ell,N}) - A^{**}(\th_{\ell-1,N})}{\th_{\ell,N} - \th_{\ell-1,N}} \\
		& = \frac{A^{**}(m'') - A^{**}(m_L'')}{m'' - m_L''} - \frac{A^{**}(m_L'') - A^{**}(m')}{m_L'' - m'} =:2\s>0.
		\end{split} 
		\]
		Thus Lemma \ref{l:subdisc*} applies.  From this point onward, the argument is identical to that of Case 1. 
	\end{proof}

	\section{Finite-time clustering in the supercritical region}
	
	\label{sec:super}
	
	In this section, we turn to the \emph{supercritical region} $\Sigma_-$, where the flux $A$ is detached from its convex envelope $A^{**}$.  We prove part II of Theorem \ref{t:regions}, which guarantees that any compact subinterval of $\Sigma_-$ becomes part of a finite-time cluster; this phenomenon is completely new in the study of \eqref{e:EA}.  We also show that clusters inside a given connected component $(m_-, m_+)$ of $\Sigma_-$ cannot protrude from $(m_-, m_+]$ unless they contain all of $(m_-, m_+]$.  This will be useful in the proof of part III(i) of Theorem \ref{t:regions}, which we give in Section \ref{sec:phidependent} below.  We begin by establishing discrete versions of these facts before bootstrapping to the full continuum versions of the statements.  
	
	As before, we fix a solution $(\rho, u)$ of \eqref{e:EA} and assume (A1)--(A4) hold.  
	
	\subsection{The discrete supercritical setting}
	
	\begin{PROP}\label{prop:collision}
		Let $(m_{i,N}, x_{i,N}^0, v_{i,N}^0)_{i=1}^N$ be a discretization of $(\rho^0, u^0)$ satisfying \eqref{eq:D2}--\eqref{e:psii0Nvi0N}.  Assume that $(\th_{i,N},\th_{j,N})$ is a connected component of $\Sigma_-$, so that $A >A^{**}$ on $(\th_{i,N},\th_{j,N})$ and 
		$A(\th_{\ell,N})=A^{**}(\th_{\ell,N})$ for $\ell = i,j$.
		Then the following statements hold:
		\begin{itemize}
			\item[(i)] Agents $i+1,\ldots, j$ will collide in finite time.
			\item[(ii)] Let $T_N$ be the first time when agents $i+1$ and $j$ collide.  Then agents $i$ and $i+1$ cannot collide before time $T_N$; similarly, agents $j$ and $j+1$ cannot collide before time $T_N$.
		\end{itemize}
	\end{PROP}
	\begin{proof}
		We start with statement (ii), which is actually a special case of something we have already proven.  If agents $i$ and $i+1$ ever collide (say at time $t\ge 0$), then the hypotheses of Lemma \ref{l:inductionlemma} are satisfied, so that in particular we must have $\th_{i^*(t),N}\notin \Sigma_-$.  As $j$ is the smallest index exceeding $i$ such that $\th_{j,N}\notin \Sigma_-$, and since $i^*(t)> i$ by assumption, we must have $i^*(t)\ge j$.  Furthermore, $T_N=\inf\{s:i^*(s)\ge j\}$, and therefore $t\ge T_N$, as claimed.  The situation is the same for agents $j$ and $j+1$.
		
		Now we prove statement (i), arguing by contradiction.  Assume that agents $i+1$ and $j$ never collide.  Then for any $t\ge 0$, we have  $(i+1)^*(t)<j$; furthermore, by statement (ii), agents $i$ and $i+1$ cannot collide, so we also have $(i+1)_*(t)=i+1$. Set 
		\[
		h_N=\min_{i+1\leq \ell\leq j-1}(A-A^{**})(\th_{\ell,N});
		\qquad 
		c_N=\frac{h_N}{\th_{j,N}-\th_{i,N}}.
		\]
		Then  
		\begin{align*}
		\psi_{i+1,N}(t)
		& = \frac{A(\th_{(i+1)^*(t),N})-A(\th_{i,N})}{\th_{(i+1)^*(t),N}-\th_{i,N}} \\
		& = \frac{(A-A^{**})(\th_{(i+1)^*(t),N}) - (A-A^{**})(\th_{i,N}) }{\th_{(i+1)^*(t),N}-\th_{i,N}} + 
		\frac{A^{**}(\th_{(i+1)^*(t),N})-A^{**}(\th_{i,N})}{\th_{(i+1)^*(t),N}-\th_{i,N}} \\
		& \ge c_N + \frac{A^{**}(\th_{j,N})-A^{**}(\th_{i,N})}{\th_{j,N}-\th_{i,N}}.
		\end{align*}
		(In the last line, we used linearity of $A^{**}$ on $[\th_{i,N}, \th_{j,N}]$).
		Similarly, 
		\[
		\psi_{j,N}(t)
		\le c_N +  \frac{A^{**}(\th_{j,N})-A^{**}(\th_{i,N})}{\th_{j,N}-\th_{i,N}}.
		\]
		Putting the two estimates above together, we obtain a strictly positive lower bound on the difference:
		\begin{equation}\label{eq:psidiff}
		\psi_{i+1,N}(t)-\psi_{j,N}(t)\geq 2 c_N>0,
		\end{equation}
		Note in particular that $c_N$ is independent of $t$.  Therefore, we can compute using \eqref{e:xdiff} that for all $t\ge 0$,
		\begin{align}
		\frac{\dd}{\dt}\big(x_{j,N}(t)-x_{i+1,N}(t)\big)=&\,\psi_{j,N}(t)-\psi_{i+1,N}(t)-\sum_{\ell=1}^N m_{\ell,N} \int_{x_{i+1,N}(t)}^{x_{j,N}(t)} \phi(y-x_{\ell,N}(t))\dy 
		\label{eq:xdiff4} \leq -2c_N <0.\nonumber
		\end{align}
		Here, we have used \eqref{eq:psidiff} and the fact that $\phi$ is nonnegative. Therefore, agents $i+1$ and $j$ must collide at some finite time $T_N\leq D^0/(2c_N)$, where $D^0 = \diam \supp(\rho^0)$. This finishes the proof.
	\end{proof}
	
	\subsection{Finite-time clustering for the Euler-alignment system}
	
	Next, we use the sticky particle approximation to extend our result to the continuum system. The difficulty is the dependence of the constant $c_N$ in \eqref{eq:psidiff} on $N$. In fact, 
	\[
	\lim_{N\to\infty}h_N=0,\quad\lim_{N\to\infty}c_N=0.
	\]
	Consequently, our upper bound on $T_N$ tends to infinity with $N$.  The argument of the previous subsection therefore does \textit{not} prove that a connected component $(m_-, m_+)$ of $\Sigma_-$ (or more precisely, its image under $X(\cdot, t)$) collapses to a point in finite time, and in fact, this is not true in general.  However, as we will show below, any compact subset $K$ of $(m_-,m_+)$ \textit{will} collapse to a point in a finite time $T$ that we can bound from above.  In order to streamline the proof of this theorem, we first state an elementary but technical lemma that will help us obtain a uniform-in-$N$ adaptation of the statement and argument of Proposition \ref{prop:collision}.
	
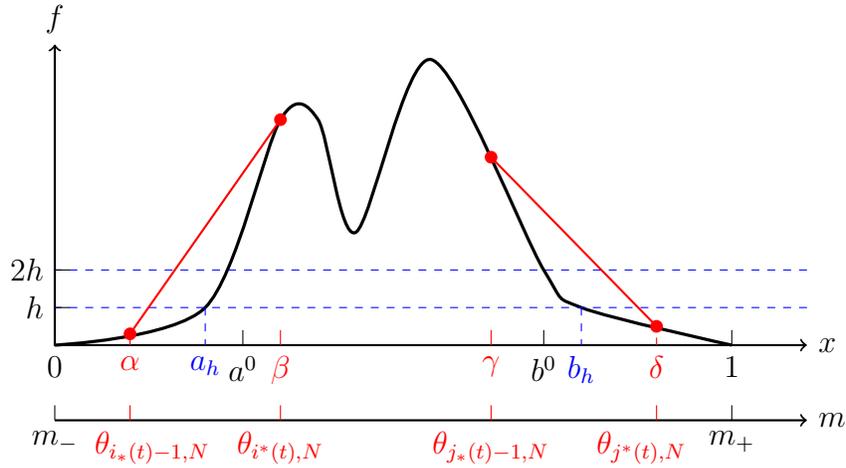
\begin{figure}[ht]
		\begin{tikzpicture}
		\draw[thick,<->] (10,0) node[right]{$x$} -- (0,0) node[below]{$0$} -- (0,4) node[above]{$f$};
		\draw[very thick] plot[smooth] coordinates {(0,0)   (2,.5) (3,3)  (3.5,3) (4,1.5) (5,3.8) (6.5, 1) (7,.5) (9,0)};

		\draw[blue, dashed] (10,.5) -- (0,.5);
		\draw (.2,.5) -- (0,.5)  node[left] {$h$};
		
		\draw[thick, red] plot[mark=*] coordinates {(1,.15) (3,3)};
		\draw[thick, red] plot[mark=*] coordinates {(8,.25) (5.8,2.5)};
		
		\draw[blue, dashed] (10,1) -- (0,1);
		\draw (.2,1) -- (0,1)  node[left] {$2h$};
		
		\draw[red] (1,.1) -- (1,0) node[below] {$\alpha$};
		\draw[blue, dashed] (2,.5) -- (2,0) node[below] {$a_h$};
		\draw (2.5,.2) -- (2.5,0) node[below]{$a^0$};
		\draw[red] (3,.2) -- (3,0) node[below]{$\beta$};
		\draw[red] (5.8,.2) -- (5.8,0) node[below]{$\gamma$};
		\draw (6.5,.2) -- (6.5,0) node[below]{$b^0$};
		\draw[blue, dashed] (7,.5) -- (7,0) node[below]{$b_h$};
		\draw[red] (8,.1) -- (8,0) node[below]{$\delta$};
		\draw (9,.2) -- (9,0) node[below]{$1$};
		
		\draw[thick,->] (0,-1) -- (10,-1) node[right]{$m$};
		\draw (0,-.8) -- (0, -1) node[below]{$m_-$};
		\draw[red] (1,-.8) -- (1,-1) node[below,xshift=.3cm] {$\th_{i_*(t)-1,N}$};
		\draw[red] (3,-.8) -- (3,-1) node[below]{$\th_{i^*(t),N}$};
		\draw[red] (5.8,-.8) -- (5.8,-1) node[below]{$\th_{j_*(t)-1,N}$};
		\draw[red] (8,-.8) -- (8,-1) node[below,xshift=-.2cm]{$\th_{j^*(t),N}$};
		\draw (9,-.8) -- (9,-1) node[below]{$m_+$};
		
		\end{tikzpicture}
		\caption{Illustration of Lemma \ref{lem:Fluxbd} (with $h=h^0$) and its (rescaled) application in the proof of Theorem \ref{thm:sup}}
		\label{fig:f}
	\end{figure}

	\begin{LEMMA}\label{lem:Fluxbd}
		Let $f:[0,1]\to \R$ be a Lipschitz function, and let $\widetilde{K}$ be a compact subset of $(0,1)$.  Assume that 
		\begin{itemize}
			\item [(a)] $f(0) = f(1) = 0$, and $f>0$ on $(0,1)$,
			\item [(b)] $f$ is convex in a neighborhood of 0, and a neighborhood of 1.
		\end{itemize}	
		Then there exists $h^0>0$ such that whenever $0<h\le h^0$, the following holds:
		\begin{itemize}
			\item [(i)] $f^{-1}(h)$ consists of exactly two points, call them $a_h$ and $b_h$, with $a_h<b_h$, and $\widetilde{K}\subseteq [a_h,b_h]$.
			\item [(ii)] $f>h$ on $(a_h,b_h)$ and $f<h$ on $[0,a_h)\cup (b_h,1]$.
			\item [(iii)] Whenever  
			\begin{equation} 
			\label{e:lemmaorder}
			0\leq\alpha<a_h\leq\beta<\gamma<b_h\leq\delta\leq1,
			\end{equation} 
			we have the following estimate
			\begin{equation}\label{eq:lemlowerbound}
			\frac{f(\beta) - f(\alpha)}{\beta - \alpha} - \frac{f(\delta) - f(\gamma)}{\delta - \gamma}\geq h>0,
			\end{equation}
		\end{itemize}
	\end{LEMMA}
	See Figure \ref{fig:f} for an illustration of the statement of the Lemma.  Later on, we will apply the Lemma with $f$ equal to a rescaled version of $A - A^{**}$.

	\begin{proof}
		First, we choose $a^0$ and $b^0$ such that $\widetilde{K}\subseteq[a^0,b^0]\subset(0,1)$, and $f$ is convex in $[0,a^0]$ and $[b^0,1]$. By assumption (a), $f$ must be strictly increasing in $[0,a^0]$ and strictly decreasing in $[b^0,1]$. 
		Since $f>0$ in $[a^0,b^0]$, its minimum on this interval must be positive. We denote 
		\begin{equation}\label{eq:h0}
		h^0=\frac12\min_{x\in[a^0,b^0]}f(x)>0.
		\end{equation}
		For any $h$ in $(0,h^0]$, we define $a_h\in(0,a^0)$ by $a_h=f^{-1}(h)$. Note that $a_h$ is well-defined, as $f$ is strictly increasing on $[0,a^0]$ and $0=f(0)<h< f(a^0)$. Similarly, we define $b_h\in(b^0,1)$ by $b_h=f^{-1}(h)$. It is clear that (i) and (ii) hold, and that furthermore, for $\a,\b,\g,\d$ satisfying \eqref{e:lemmaorder}, we must have
		\[
		\frac{f(\beta) - f(\alpha)}{\beta - \alpha} > 0 > \frac{f(\delta) - f(\gamma)}{\delta - \gamma},
		\qquad \text{ which implies } 
		\qquad 
		\frac{f(\beta) - f(\alpha)}{\beta - \alpha} - \frac{f(\delta) - f(\gamma)}{\delta - \gamma}>0.
		\]
		We need to improve this to a positive lower bound that is uniform in $\a,\b,\g,\d$.  We do this by considering the following three cases.
		
		\noindent\textit{Case 1: $\beta<a^0$}.  Using the convexity of $f$ on $[0,a^0]$, we have
		\[\frac{f(\beta) - f(\alpha)}{\beta - \alpha} - \frac{f(\delta) - f(\gamma)}{\delta - \gamma}\geq\frac{f(a) - f(0)}{a - 0}-0=\frac{h}{a}>h.\]
		\noindent\textit{Case 2: $a^0\leq\beta<\gamma\leq b^0$}. From the definition of $h^0$ in \eqref{eq:h0}, we have $f(\beta)\ge 2h_0$ and $f(\gamma)\ge 2h^0$. Hence,
		\[
		\frac{f(\beta) - f(\alpha)}{\beta - \alpha} - \frac{f(\delta) - f(\gamma)}{\delta - \gamma}\geq\frac{2h^0 - h}{1 - 0}-\frac{h-2h^0}{1-0}\ge 2h.
		\]
		\noindent\textit{Case 3: $\gamma>b^0$}. Using the convexity of $f$ on $[b^0,1]$, we obtain 
		\[\frac{f(\beta) - f(\alpha)}{\beta - \alpha} - \frac{f(\delta) - f(\gamma)}{\delta - \gamma}\geq0-\frac{f(1) - f(b)}{1 - b}=\frac{h}{1-b}>h.\]
		This completes the proof.
	\end{proof}
	\begin{REMARK}
		The technical assumption (b) on the convexity near the boundary is to make sure that there exists a small $h$ such that (ii) holds.  We want to eliminate the possibility of highly oscillatory functions like 
		\[f(x) = x(2+\sin\tfrac{1}{x})\]
		near $0$.  We will not attempt to treat fluxes that exhibit this sort of pathological behavior.
	\end{REMARK}
	
	We now prove part II of Theorem \ref{t:regions}, restated here in slightly modified but equivalent form.
	\begin{THEOREM}
		\label{thm:sup}
		Let $(\rho, u)$ and the accompanying functions and sets be as in (A1)--(A4).  Assume $(m_-,m_+)$ is a connected component of $\Sigma_-$, i.e., 
		\[
		A(m)>A^{**}(m),\quad\forall~m\in (m_-,m_+),\quad\text{and}\quad
		A(m_-)=A^{**}(m_-),\quad A(m_+)=A^{**}(m_+).
		\]
		We also recall that (A4) requires $A$ to be convex in a neighborhood of $m_-$ and $m_+$.  Let $K$ be a compact subset of $(m_-, m_+)$. Then there exists a finite time $T>0$ such that 
		\[
		X(m,t) = X(m',t)\quad\forall~m,m'\in K\quad\text{and}\quad \forall~t\geq T.
		\]
	\end{THEOREM}
	
	\begin{proof}
		Our plan is to approximate the system by the discrete dynamics, apply Lemma \ref{lem:Fluxbd} to obtain uniform bounds, and then pass to the limit.
		
		\textbf{Step 1:} Define the linear bijection  $m:[0,1]\to [m_-,m_+]$ via 
		\[
		m(x) = m_- + (m_+ - m_-)x.
		\]
		Then define $f:[0,1]\to\R$ by
		\begin{equation}\label{def:f}
		f(x)=A(m(x))-A^{**}(m(x)),
		\end{equation}
		and define $\widetilde{K}\subset (0,1)$ by $\widetilde{K} = m^{-1}(K)$.  
		It is easy to check that $f$ satisfies all the assumptions in Lemma~\ref{lem:Fluxbd}, as a consequence of assumption (A4). Choose $h^0$ as in the conclusion of the Lemma.  Then for any $h\in (0,h^0]$, we have that 
		\[
		K\subseteq [m(a_h), m(b_h)],
		\]
		and whenever we have 
		\begin{equation*}
		m_- < \a < m(a_h) \le \b < \g < m(b_h) \le \d \le m_+,
		\end{equation*}
		then the following estimate holds:
		\begin{equation}
		\label{e:lemmaapplied}
		\frac{(A-A^{**})(\b) - (A-A^{**})(\a)}{\b - \a} -  \frac{(A-A^{**})(\d) - (A-A^{**})(\g)}{\d - \g} \ge \frac{h}{m_+ - m_-}=:c(h)>0.
		\end{equation}
		Simplifying the left side of \eqref{e:lemmaapplied} using the linearity of $A^{**}$ on $[m_-, m_+]$, we obtain
		\begin{equation}
		\label{e:lemmaapplied2}
		\frac{A(\b) - A(\a)}{\b - \a} -  \frac{A(\d) - A(\g)}{\d - \g} \ge c(h)>0.
		\end{equation}

		
		\textbf{Step 2:} Form a sequence of approximating systems $(m_{i,N},x_{i,N}(\cdot),v_{i,N}(\cdot))_{i=1}^N$ via the procedure in Section \ref{sec:preliminaries}. (In particular, assume \eqref{eq:milimit}--\eqref{e:psii0Nvi0N} are satisfied).  For $N$ large enough, assumption \eqref{eq:D2} says that $S\subset \{\th_{i,N}\}_{i=0}^N$, so that in particular $m_-$ and $m_+$ are breakpoints.  
		
		Define 
		\begin{equation}
		\label{e:collisiontimebd}
		T:=\frac{X^0(m_+) - X^0(m_-)}{c(\frac{h_0}{2})}, 
		\end{equation}
		where $c(h) = \frac{h}{m_+ - m_-}$, as defined in \eqref{e:lemmaapplied}. Fix any $t\ge T$, and choose a subsequence $(N_k)_{k=1}^\infty$ such that the a.e. $\!\!$ convergence \eqref{e:XNkconv} holds at time $t$.  Then we may in particular choose $h\in (\frac{h^0}{2},h^0]$ so that, putting 
		\begin{equation}
		\label{e:thetah}
		\th' = m(a_h), 
		\qquad \th'' = m(b_h),
		\end{equation}
		we have
		\begin{equation}\label{eq:XNnconv4}
		\begin{split}
		X_{N_k}(\th',t)\to X(\th',t)
		\qquad \text{ and } \qquad 
		X_{N_k}(\th'',t)\to X(\th'',t), 
		\end{split}
		\qquad \text{ as } k\to \infty.
		\end{equation}
		
		Note that the subsequence $(N_k)_{k=1}^\infty$, and hence our choice of $h$, depends on $t$.  However, $h^0$ is of course time-independent.
		
		Our goal is to show that  
		\begin{equation}\label{eq:Xcluster4}
		X(\th',t)=X(\th'',t).
		\end{equation}
		Since $K\subseteq [\th', \th'']$ and $t\ge T$ is arbitrary, establishing \eqref{eq:Xcluster4} will prove the theorem.
		
		\textbf{Step 3:} For each $N\in \N$, choose $i,j\in \{1, \ldots, N\}$ such that $\th'\in (\th_{i-1,N},\th_{i,N}]$ and $\th''\in (\th_{j-1,N},\th_{j,N}]$.  
		Using \eqref{e:thetah} and \eqref{e:lemmaapplied2} (with $\a = \th_{i_*(s)-1,N}$, $\b = \th_{i^*(s),N}$, $\g = \th_{j_*(s)-1,N}$, and $\d = \th_{j^*(s),N}$), we have the following uniform lower bound for all $s\ge 0$:
		\begin{align*}
		\psi_{i,N}(s)-\psi_{j,N}(s)=&\,\frac{A(\th_{i^*(s),N}) -A(\th_{i_*(s)-1,N})}{\th_{i^*(s),N}-\th_{i_*(s)-1,N}} -\frac{A(\th_{j^*(s),N})-A(\th_{j_*(s)-1,N})}{\th_{j^*(s),N}-\th_{j_*(s)-1,N}}
		\ge c(h)\ge c(\tfrac{h_0}{2})>0.
		\end{align*}
		It follows that for each $N\in \N$ and each time $s$ prior to the collision of agents $i$ and $j$, we have
		\begin{equation}\label{eq:xdiffode4}
		\frac{\dd}{\ds}(x_{j,N}(s)-x_{i,N}(s)) \leq-\big(\psi_{i,N}(s)-\psi_{j,N}(s)\big)\le -c(\tfrac{h_0}{2}).
		\end{equation}
		Integrating this equation and recalling the definition of $T$, as well as the fact that $X_N(\th'', \cdot) - X_N(\th',\cdot)=x_{j,N}-x_{i,N}$, we deduce that 
		\[
		X_N(\th'',s)-X_N(\th',s) = 0,
		\qquad \forall s\ge T.
		\]
		Then, setting $s = t$ and taking $N\to \infty$ along the subsequence $(N_k)_{k=1}^\infty$ finishes the proof.
	\end{proof}
	
	\subsection{Confinement of `incomplete' clusters in supercritical intervals}
	
	We now upgrade Proposition \ref{prop:collision}(ii) to the full continuum version, as previewed at the beginning of this section.
	
	\begin{THEOREM}
		\label{t:confine}
		Suppose $(m',m'']$ is a $T$-cluster containing $m\in \Sigma_-$, and let $(m_-, m_+)$ be the connected component of $\Sigma_-$ containing $m$.  If $(m',m'']$ contains points outside of $(m_-, m_+]$ then $(m',m'']$ must also contain all of $(m_-, m_+]$.  
	\end{THEOREM}
	
	\begin{proof}	
		We assume, for purposes of obtaining a contradiction, that $m_-<m'<m_+<m''$.  (The case where $m'<m_-<m''<m_+$ is similar.)  Since $A$ is not linear on $(m',m'']$, the function $X^0$ cannot be constant on this interval, so we may assume without loss of generality that $T>0$.
		
		\textbf{Step 1.}  We define the relevant parameters as follows. 	 Choose $a^0,b^0\in (m_-, m_+)$ such that $A$ is convex on $[m_-, a^0]$ and $[b^0, m_+]$.  Put $h^0 = \frac12 \min_{[a^0,b^0]} (A - A^{**})$; for any $h\in (0,2h^0)$, let $a_h\in (m_-, a^0)$ and $b_h\in (b^0, m_+)$ denote the two points in the inverse image of $h$ under $A - A^{**}$ that lie in $(m_-, m_+)$. Take $h^1\in (0,h^0)$ small enough so that $a_{2h^1}<m'<b_{h^1}$.
		
		As usual, consider a sequence of discretized systems satisfying \eqref{eq:milimit}--\eqref{e:psii0Nvi0N}, and additionally $m_-, m_+\in \{\th_{i,N}\}_{i=1}^N$ for sufficiently large $N$.  Let $(N_k)_{k=1}^\infty$ be a subsequence such that the a.e. $\!\!$ convergence \eqref{e:XNkconv} holds at time $t = T$.  Choose $\th'\in (b_{h^1}, m_+)$ and $\th''\in (m_+, m'')$ such that 
		\[
		X_{N_k}(\th',T)\to X(\th', T) 
		\qquad \text{ and } \qquad 
		X_{N_k}(\th'',T)\to X(\th'', T),
		\qquad \text{ as } k\to \infty.
		\] 
		
		Fix $h\in (0,h^1)$ such that $b_h = \th'$.  And finally, for each sufficiently large $N\in \N$, choose indices $I<i<J$ such that
		\[
		m_- = \th_{I,N},
		\qquad 
		\th'\in (\th_{i-1,N}, \th_{i,N}]\subset (m', m_+),
		\qquad 
		m_+ = \th_{J,N}.
		\] 
		For convenience, we review the ordering of the relevant quantities before proceeding:
		\[
		m_- = \th_{I,N} < a_{2h} < m' < \th_{i-1,N} < b_h = \th' < \th_{i,N} < m_+ = \th_{J,N} < \th'' < m''.
		\]
		
		\textbf{Step 2.} We will now reason similarly to the proof of Lemma \ref{lem:Fluxbd} and prove that 
		\begin{equation}
		\label{e:A-A**}
		\frac{(A-A^{**})(\th_{i^*(t),N_k}) - (A-A^{**})(\th_{i_*(t)-1,N_k})}{\th_{i^*(t),N_k} - \th_{i_*(t)-1,N_k}} < -\frac{h}{m_+ - m_-}.	
		\end{equation}
		
		Since $\th'$ lies in the $T$-cluster $(m',m'']$ at $m$ but $a_{2h}<m'$ does not, we have $X(\th',T)>X(a_{2h},T)$; consequently, we may assume that $X_{N_k}(\th',T)>X_{N_k}(a_{2h},T)$ for all $k$, so that 
		\begin{equation}
		\label{e:a2h}
		\th_{i_*(t)-1,N_k}\ge a_{2h}
		\qquad \forall t\in [0,T], \;\forall k\in \N.
		\end{equation}
		There are two cases to consider.  If $\th_{i_*(t)-1,N_k}\le b_{2h}$, then 
		\[
		(A-A^{**})(\th_{i_*(t)-1,N_k})\ge 2h
		\qquad \text{ and } \qquad 
		(A-A^{**})(\th_{i^*(t),N_k})< h.
		\]
		The inequality \eqref{e:A-A**} follows immediately.  On the other hand, if $b_{2h}<\th_{i_*(t)-1,N_k}<b_h$, then by convexity of $A-A^{**}$ on $[b^0, m_+]$ we get 
		\[
		\frac{(A-A^{**})(\th_{i^*(t),N_k}) - (A-A^{**})(\th_{i_*(t)-1,N_k})}{\th_{i^*(t),N_k} - \th_{i_*(t)-1,N_k}}
		\le \frac{(A-A^{**})(m_+) - (A-A^{**})(b_h)}{m_+ - b_h} < -\frac{h}{m_+ - m_-}, 
		\]
		as claimed.  
		
		\textbf{Step 3.} We will use \eqref{e:A-A**}, together with Proposition \ref{prop:collision}(ii), to prove the following lower bound:
		\[
		\psi_{J+1,N}(t) - \psi_{i,N}(t) \ge \frac{h}{m_+ - m_-} =:2\s>0.
		\]
		This is, in a sense, the key step of the proof, since it is here that we leverage the discrete version of the statement we want to prove.  
		
		Recalling \eqref{e:a2h}, we know that 
		\[
		\th_{I,N_k}<a_{2h}<\th_{i_*(t)-1,N_k}\le\th_{J_*(t)-1,N_k},
		\qquad \forall t\in [0,T],
		\]
		which implies that $J_*(t)>I+1$, i.e., agents $I+1$ and $J$ do not collide on the time interval $[0,T]$.  Consequently, Proposition \ref{prop:collision}(ii) tells us that agents $J$ and $J+1$ also do not collide on $[0,T]$.  Therefore, for any time $t\in [0,T]$, we thus have the following:
		\begin{align*}
		\psi_{i,N_k}(t) 
		& = \frac{A(\th_{i^*(t),N_k}) - A(\th_{i_*(t)-1,N_k})}{\th_{i^*(t),N_k} - \th_{i_*(t)-1,N_k}}\\
		& = \frac{(A-A^{**})(\th_{i^*(t),N_k}) - (A-A^{**})(\th_{i_*(t)-1,N_k})}{\th_{i^*(t),N_k} - \th_{i_*(t)-1,N_k}} + \frac{A^{**}(\th_{i^*(t),N_k}) - A^{**}(\th_{i_*(t)-1,N_k})}{\th_{i^*(t),N_k} - \th_{i_*(t)-1,N_k}}\\
		& \le -\frac{h}{m_+ - m_-} + \frac{A^{**}(\th_{(J+1)^{*}(t),N_k}) - A^{**}(\th_{J,N_k})}{\th_{(J+1)^{*}(t),N_k} - \th_{J,N_k}} \\
		& \le -\frac{h}{m_+ - m_-} + \frac{A(\th_{(J+1)^{*}(t),N_k}) - A(\th_{J,N_k})}{\th_{(J+1)^{*}(t),N_k} - \th_{J,N_k}} =  -\frac{h}{m_+ - m_-} + \psi_{J+1,N_k}(t).
		\end{align*}
		
		This proves the desired lower bound.
		
		\textbf{Step 4.} We derive a contradiction.  Using \eqref{e:xdiff} and arguing as in (the beginning of) the proof of Lemma \ref{l:subdisc*}, we obtain 
		\[
		x_{J+1,N_k}(t) - x_{i,N_k}(t) \ge \min\{x_{J+1,N_k}^0 - x_{i,N_k}^0 + t\s, \eta\},
		\qquad \forall t\ge 0,
		\]
		where $\eta>0$ is chosen small enough so that $|\int_z^w \phi(y)\dr|<\s$ when $|z - w|<\eta$.  It follows that 
		\[
		X_{N_k}(\th'',T) - X_{N_k}(\th',T) \ge x_{J+1,N_k}(T) - x_{i,N_k}(T)\ge \min\{T\s, \eta\}>0.
		\] 
		Taking $k\to \infty$ and recalling that $m'<\th'<\th''<m''$, we conclude that 
		\[
		X(m'',T) - X(m',T)
		\ge \min\{T\s, \eta\}>0,
		\]
		which contradicts our initial assumption that $(m',m'']$ is a $T$-cluster and therefore finishes the proof of the theorem.	
	\end{proof}

	\section{Protocol-dependent behavior:  clustering and non-clustering in the critical regime}
	
	\label{sec:phidependent}
	
	So far, we have shown that $X(\cdot, t)$ exhibits two distinct behaviors in the subcritical region $\Sigma_+$ and the supercritical region $\Sigma_-$. We have also shown that clustering does not occur across different $L(m)$'s.  In this section, we study the clustering behavior within a single  $L(m)$.  In contrast to the analysis we have presented thus far, the results of this section depend on properties of $\phi$ beyond those assumed in (A1); in that sense, the situation we consider in this section constitutes a sort of `critical' regime (which includes in particular the behavior in the entire critical region $\Sigma_0$ but may also concern intervals in $\Sigma_-$).
	
	The three subsections here track statements (i), (ii), and (iii) of part III of Theorem \ref{t:regions}.  The first statement, which concerns bounded protocols $\phi$, places limitations on the kinds of finite-time clusters that can occur within a given $L(m)$ and is the most delicate of the three results discussed in this section. The reason is that, unlike the situation in Section \ref{sec:sub}, it is possible for mass labels \textit{within} a given $L(m)$ to belong to the same infinite-time cluster without belonging to the same $t$-cluster for any finite time.  Ruling out finite-time clustering is therefore more subtle here than in the subcritical regime.  The final two subsections show that under the heavy-tail assumption \eqref{eq:heavytail}, each $L(m)$ becomes an infinite-time cluster, and is actually a finite-time cluster if $\phi$ is weakly singular.  
	
	As always, we assume throughout this section that $(\rho, u)$ is an entropy solution of \eqref{e:EA} and that (A1)--(A4) hold.
	
	\subsection{Bounded communication}
	
	In this subsection, we assume that $\phi$ is bounded, and we prove statement III(i) of Theorem \ref{t:regions}.  As usual, we will start with a discrete (partial) analog of our target statement.  
	
	\subsubsection{The discrete setting} 
	
	\begin{LEMMA}
		\label{l:disccritlwrbd}
		Suppose that $\phi\in L^\infty(\R)$.  Let $(m_{i,N}, x_{i,N}^0, v_{i,N}^0)_{i=1}^N$ be a discretization of $(\rho^0, u^0)$ satisfying \eqref{eq:D2}--\eqref{e:psii0Nvi0N}.  Assume that $i,j\in \{1, \ldots, N\}$ are such that $i\le j$ and $\th_{i,N},\th_{j,N}\notin \Sigma_-$.  
		Then 
		\begin{equation}
		\label{e:disccritlwrbd}
		x_{j+1,N}(t) - x_{i,N}(t)\ge (x_{j+1,N}^0 - x_{i,N}^0)e^{-\|\phi\|_\infty t}, 
		\qquad \forall t\ge 0.
		\end{equation}	
	\end{LEMMA}
	
	\begin{proof}
		Assume without loss of generality that $x_{i,N}^0 < x_{j+1,N}^0$.  Since $\th_{i,N}, \th_{j,N}\notin \Sigma_-$, Corollary \ref{c:Apinned} guarantees that 
		\[
		\psi_{j+1,N}(t)-\psi_{i,N}(t)\ge 0,
		\qquad \forall t\ge 0.
		\]  
		
		Since $\phi$ is bounded, we can use an improved estimate on the sum in \eqref{e:xdiff}:
		\begin{align*}
		\sum_{\ell=1}^N m_{\ell,N}\int_{x_{i,N}(s)}^{x_{j+1,N}(s)}\phi(y-x_{\ell,N}(s))\dy
		\le \|\phi\|_{\infty}\big(x_{j+1,N}(s)-x_{i,N}(s)\big).
		\end{align*}
		Consequently,
		\[	
		\frac{\dd}{\ds}\big(x_{j+1,N}(s)-x_{i,N}(s)\big)\geq- \|\phi\|_{\infty}\big(x_{j+1,N}(s)-x_{i,N}(s)\big),
		\]
		which becomes \eqref{e:disccritlwrbd}  after integration.  
	\end{proof}
	
	\subsubsection{The continuum setting}
	
	We now turn to the proof of statement III(i) in Theorem \ref{t:regions}.  Under the assumption that $\phi\in L^\infty(\R)$, we want to show that there are only two possible kinds of clusters.  First, there might be initial clusters---intervals $(m',m'']$ on which $X^0$ is constant.  And second, there might be clusters that are contained in some $(m_-,m_+]$, where $(m_-,m_+)$ is a connected component of $\Sigma_-$.  We will show that no other finite-time clusters are possible.  With this in mind, we make the following notation, which we will use in the proof of our theorem.
	\begin{DEF}
		Fix $m\in (-\frac12, \frac12]$.  Define the set $C(m)$ as follows.  
		\begin{itemize}
			\item If there is an initial cluster at $m$, let $C(m)$ denote this initial cluster.
			\item If $m\in \Sigma_-$, define $C(m) = (m_-, m_+]$, where $(m_-, m_+)$ is the connected component of $\Sigma_-$ that contains $m$.
			\item Otherwise, put $C(m) = \{m\}$.
		\end{itemize}	
		We also define the function $R(m,t)$ as follows.
		\[
		R(m,t) = \begin{cases}
		\frac{1}{|C(m)|}\displaystyle\int_{C(m)} X(\widetilde{m},t)\,\dd\widetilde{m},
		& \text{ if } m\in \Sigma_- \\
		X(m,t), & \text{ otherwise.} \\
		\end{cases}
		\]
		Here $|C(m)|$ denotes the Lebesgue measure of the interval $C(m)$.  
	\end{DEF}
	
	The following theorem is a more precise version of statement III(i) in Theorem \ref{t:regions}. 
	
	\begin{THEOREM}\label{thm:critbound1}
		Assume $\phi\in L^\infty(\R)$.
		\begin{itemize}
			\item [(i)] If $C(m'') = \{m''\}$ is a singleton and $m'<m''$, then there exists a constant $c>0$ such that 
			\begin{equation}
			\label{e:Xsepcrit}
			X(m'',t) - X(m',t)\ge c e^{-\|\phi\|_{\infty} t},
			\qquad \forall t\ge 0.
			\end{equation}
			\item [(ii)] If $C(m'')$ is an interval and $m'<\inf C(m'')$, then 
			\begin{equation}
			\label{e:Rsep}
			R(m'',t) - R(m',t)\ge c e^{-\|\phi\|_{\infty} t},
			\qquad \forall t\ge 0.
			\end{equation}	
		\end{itemize} 
		In either case, $m'$ and $m''$ do not belong to the same $t$-cluster for any finite time $t$.
	\end{THEOREM}

	\begin{proof}
		Note first of all that Theorem \ref{t:confine} is the link between \eqref{e:Rsep} and the last claim, in the event that $C(m'')$ is an interval: If $m'$ and $m''$ ever belonged to the same $t$-cluster, then all of $C(m')$ and $C(m'')$ would need to belong to that $t$-cluster as well, by Theorem \ref{t:confine}.  This would in turn imply that $R(m',t) = R(m'',t)$, in direct contradiction with~\eqref{e:Rsep}.	
		
		We also note that, by Theorem \ref{thm:sub}, we may assume without loss of generality that $m'\in L(m'')$.  
		
		\textbf{Proof of (i).} Assume that $C(m'') = \{m''\}$. Then we must have $m''\in \Sigma_0$ (since $m''\in \Sigma_+$ would force $L(m')\ne L(m'')$), so we may assume without loss of generality that the interval $[m',m'']$ lies entirely in $\Sigma_0$.  Since $X^0$ is left-continuous and $X^0$ is not constant on any interval of the form $[m'' - \e, m'']$ for $\e>0$, it follows that $X^0$ must take infinitely many values on the interval $[m',m_-'']$.  We may therefore choose $\widetilde{m}',\widetilde{m}''\in [m',m'']$ such that 
		\[
		X^0(m')\le X^0(\widetilde{m}') < X^0(\widetilde{m}'')\le X^0(m'')
		\]
		Consider as usual a sequence of discretizations satisfying \eqref{eq:milimit}--\eqref{e:psii0Nvi0N}; fix a time $t>0$, and choose a corresponding subsequence $(N_k)_{k=1}^\infty$ such that the a.e. $\!\!$ convergence \eqref{e:XNkconv} holds at time $t$.  Choose $\th'\in (m', \widetilde{m}')$, $\th''\in (\widetilde{m}'', m_-'')$ such that 
		\begin{equation}
		\label{e:XNKCmc1}
		\begin{split} 
		X_{N_k}(\th',t)\to X(\th',t),
		\qquad 
		& X_{N_k}(\th'',t)\to X(\th'',t), \\
		X_{N_k}^0(\th')\to X^0(\th'),
		\qquad 
		& X_{N_k}^0(\th'')\to X^0(\th''),
		\end{split} 
		\qquad \text{ as } k\to \infty.
		\end{equation}
		For each sufficiently large $N\in \N$, choose $i,j\in \{1, \ldots, N\}$ such that 
		\[
		\th'\in (\th_{i-1,N}, \th_{i,N}]\subset (m',\widetilde{m}');
		\qquad 
		\th''\in (\th_{j,N}, \th_{j+1,N}]\subset (\widetilde{m}'',m_-'').
		\]
		Since $[m',m'']\subset \Sigma_0$ (and in particular $[m',m'']$ does not intersect $\Sigma_-$), Lemma \ref{l:disccritlwrbd} implies that 
		\begin{equation}
		x_{j+1,N}(t) - x_{i,N}(t)\ge (x_{j+1,N}^0 - x_{i,N}^0) e^{-\|\phi\|_{\infty} t}
		\end{equation} 
		Recalling that $x_{i,N} = X_N(\th',\cdot)$ and $x_{j+1,N} = X_N(\th'', \cdot)$, we deduce that
		\[
		X_{N}(\th'', t) - X_{N}(\th', t) \geq \big(X_{N}^0(\th'') - X_{N}^0(\th')\big)e^{-\|\phi\|_{L^\infty} t}.
		\]
		Then, taking $N\to\infty$ along the subsequence $(N_k)_{k=1}^\infty$ and recalling that 
		\[
		m'<\th'<\widetilde{m}'<\widetilde{m}''<\th''<m''\le m'',
		\] 
		we obtain the desired estimate, with $c = X^0(\widetilde{m}'') - X^0(\widetilde{m}')>0$.
		
		\textbf{Proof of (ii).}  Assume that $C(m'')$ is an interval and that $m'<\inf C(m'')$.  
		Define
		\[
		m_-'' = \inf C(m'').
		\]
		If $C(m_-'')$ is a singleton, i.e.,  $C(m_-'') = \{m_-''\}$,  then $m'<\inf C(m'')=m_-''$ is upgraded to $m_+':=\sup C(m')<m_-''$.  Then statement (i) can be applied, with $m_-''$ replacing $m''$, and $m_+'$ replacing $m'$, to yield
		\[
		R(m'',t) - R(m',t)\ge X(m_-'',t) - X(m_+',t)\ge ce^{-\|\phi\|_\infty t},
		\qquad \forall t\ge 0.
		\]
		Thus, we may assume that $C(m_-'')$ and $C(m'')$ are both intervals.  In this setting, we may further assume that $m'$ lies in the interval $C(m_-'')$, which implies that $C(m')$ and $C(m'')$ are adjacent half-open intervals.  Let us therefore write 
		\[
		m_-' = \inf C(m'), 
		\qquad m_+' = \sup C(m') = \inf C(m'') = m_-'',
		\qquad 
		m_+'' = \sup C(m'').
		\]
		Consider a discretization satisfying \eqref{eq:milimit}--\eqref{e:psii0Nvi0N}, which additionally satisfies (for large enough $N$, and $i,j,\ell$ depending on $N$)
		\begin{equation}
		\th_{i,N} = m_-', 
		\qquad 
		\th_{j,N} = m_+' = m_-'',
		\qquad 
		\th_{\ell,N} = m_+''.
		\end{equation}
		Note that $\th_{i,N}$, $\th_{j,N}$, and $\th_{\ell,N}$ all belong to the complement of $\Sigma_-$, so $A(\th_{p,N}) = A^{**}(\th_{p,N})$ for $p = i,j,\ell$.  Since $A^{**}$ is linear on $[\th_{i,N}, \th_{\ell,N}]$, it follows that 
		\[
		\frac{A(\th_{j,N}) - A(\th_{i,N})}{\th_{j,N} - \th_{i,N}}
		= \frac{A(\th_{\ell,N}) - A(\th_{j,N})}{\th_{\ell,N} - \th_{j,N}}.
		\]
		Furthermore, Lemma \ref{l:disccritlwrbd} implies that no collisions occur between agents $i$ and $i+1$, $j$ and $j+1$, or $\ell$ and $\ell+1$.  Therefore, the following equality holds for all time:
		\begin{equation}
		\label{e:psisumcrit}
		\frac{\sum_{p=i+1}^j m_{p,N} \psi_{p,N}(t)}{\sum_{p=i+1}^j m_{p,N}} 
		= \frac{A(\th_{j,N}) - A(\th_{i,N})}{\th_{j,N} - \th_{i,N}} 
		= \frac{A(\th_{\ell,N}) - A(\th_{j,N})}{\th_{\ell,N} - \th_{j,N}} 
		= \frac{\sum_{q=j+1}^\ell m_{q,N} \psi_{q,N}(t)}{\sum_{q=j+1}^\ell m_{q,N}}.
		\end{equation}
		The quantities on the left and right side of \eqref{e:psisumcrit} are analogous to $\psi_i$ and $\psi_j$ in the inequality \eqref{e:xdiff}.  The analogs of $x_i$ and $x_j$ are 
		\begin{equation}
		\begin{split} 
		R_N(m',t)
		& :=\frac{1}{|C(m')|}\int_{C(m')} X_N(\widetilde{m},t)\,\dd\widetilde{m}
		= \frac{\sum_{p=i+1}^j m_{p,N} x_{p,N}(t)}{\sum_{p=i+1}^j m_{p,N}}, \\
		R_N(m'',t) 
		& :=\frac{1}{|C(m'')|}\int_{C(m'')} X_N(\widetilde{m},t)\,\dd\widetilde{m}
		= \frac{\sum_{q=j+1}^\ell m_{q,N} x_{q,N}(t)}{\sum_{q=j+1}^\ell m_{q,N}}.
		\end{split}
		\end{equation}
		
		Using equation \eqref{e:psisumcrit} and performing some elementary manipulations gives us the identity
		\[
		\frac{\dd}{\dt}\big( R_N(m'',t) - R_N(m',t) \big) 
		= - \frac{\sum_{r = 1}^N \sum_{p=i+1}^j \sum_{q = j+1}^\ell m_{p,N} m_{q,N} m_{r,N} \int_{x_{p,N}(t)}^{x_{q,N}(t)} \phi(y - x_{r,N}(t))\dy}{\sum_{p=i+1}^j \sum_{q = j+1}^\ell m_{p,N}}
		\]
		Using the $L^\infty$ bound on $\phi$ and integrating then yields 
		\[
		R_N(m'',t) - R_N(m',t) \ge \big( R_N(m'',0) - R_N(m',0) \big) e^{-\|\phi\|_\infty t},
		\qquad t\ge 0.
		\]
		By the $L^1$-convergence of $(X_N(s))_{N=1}^\infty$ to $X(s)$ at $s = t$ and $s = 0$, we obtain \eqref{e:Rsep}, with $c = R(m'',0) - R(m',0)$.  	
	\end{proof}

	\subsection{Heavy-tailed communication: infinite-time clustering}
	
	The lower bounds in the previous theorem are strictly positive for any finite time $t$, indicating an absence of finite-time cluster formation within~$\Sigma_0$.  However, the theorem is silent with regard to the possibility of infinite-time cluster formation. The infinite-time clustering phenomenon was studied extensively in \cite{lear2022geometric} for the case where $\Sigma_-=\emptyset$ and the velocity field is at least $C^1$. Here we present a generalization for any entropic solution.  Just as in \cite{lear2022geometric}, we assume that $\phi$ satisfies the heavy-tail condition \eqref{eq:heavytail}.  This assumption is  sufficient to guarantee uniform global communication for all time:
	\begin{equation}\label{eq:phimin}
	\phi(x-y)\geq\underline{\phi}>0,\quad\forall~x,y\in \supp(\rho(\cdot,t)),\quad \forall~t\geq0.	
	\end{equation}
	Indeed, when \eqref{eq:heavytail} holds, it is shown in \cite[Theorem 7.2]{leslie2023sticky} that the entropic solution to the Euler-alignment system \eqref{e:EA} experiences \emph{flocking}.  That is, there exists a time-independent constant $\overline{D}$ such that
	\[ 
	\diam \supp(\rho(\cdot,t)) \le \overline{D},\quad\forall~t\ge 0,	
	\]
	The assertion \eqref{eq:phimin} follows immediately, with $\underline{\phi}=\phi(\overline{D})>0$.
	
	
	\begin{THEOREM}\label{thm:critbound2}
		Assume that $\phi$ is heavy-tailed, and suppose that $L(m') = L(m'')$, with $m'<m''$.  Then 
		\begin{equation}
		\label{e:asymptoticcluster}
		X(m'',t) - X(m', t) \le D^0 e^{-\uphi t}
		\qquad \forall t\ge 0,
		\end{equation}
		where $D^0 = \diam \supp \rho^0$ and $\uphi$ is defined as in \eqref{eq:phimin}.  Consequently, $L(m') = L(m'')$ is an infinite-time cluster.
	\end{THEOREM}
	
	\begin{proof}\textbf{Step 1:} Denote 
		\[
		m_- = \inf L(m'), 
		\qquad m_+ = \sup L(m').
		\]
		We note first of all that it suffices to prove the bound \eqref{e:asymptoticcluster} in the case $m''< m_+$, by left-continuity of $X(\cdot, t)$.  
		
		Consider a sequence of discretizations satisfying \eqref{eq:milimit}--\eqref{e:psii0Nvi0N}, and additionally, $m_-, m_+\in \{\th_{i,N}\}_{i=0}^N$.  For each $N$, choose $i,j\in \{0,\ldots, N\}$ such that $\th_{i,N}=m_- \notin \Sigma_-$ and $\th_{j,N} = m_+ \notin \Sigma_-$.  Then using  Corollary~\ref{c:Apinned} and the linearity of $A^{**}$ on $L(m')$, we conclude that
		\begin{align*}
		\psi_{i+1,N}(s) 
		& \ge \frac{A^{**}(\th_{i+1,N}) - A^{**}(\th_{i,N})}{\th_{i+1,N} - \th_{i,N}}
		= \frac{A^{**}(\th_{j,N}) - A^{**}(\th_{j-1,N})}{\th_{j,N} - \th_{j-1,N}} \ge \psi_{j,N}(s),
		\qquad \forall s \ge 0.
		\end{align*}
		Next, we use the estimate  \eqref{eq:phimin} to obtain
		\begin{equation}\label{eq:phiconvrho5}
		\sum_{\ell=1}^N m_{\ell,N}\int_{x_{i+1,N}(s)}^{x_{j,N}(s)}\phi(y-x_{\ell,N}(s))\dy
		\geq 
		\underline{\phi}\big(x_{j,N}(s)-x_{i+1,N}(s)\big).	
		\end{equation}
		Therefore, we have
		\[
		\frac{\dd}{\ds}\big(x_{j,N}(s)-x_{i+1,N}(s)\big)\leq-\underline{\phi}\big(x_{j,N}(s)-x_{i+1,N}(s)\big),
		\]
		whence
		\begin{equation}\label{eq:xdiff5}
		x_{j,N}(s)-x_{i+1,N}(s)\leq \big(x_{j,N}^0-x_{i+1,N}^0\big)e^{-\underline{\phi}t}\leq D^0e^{-\underline{\phi}s},
		\qquad \forall s\ge 0.
		\end{equation}

		\textbf{Step 2:} Fix a time $t>0$ and a corresponding subsequence $(N_k)_{k=1}^\infty$ such that the a.e. $\!\!$ convergence \eqref{e:XNkconv} holds at time $t$.  Choose $\th\in (m_-, m')$ and $\th''\in (m'',m_+)$ such that 
		\[
		X_{N_k}(\th', t)\to X(\th', t) 
		\qquad \text{ and }  \qquad 
		X_{N_k}(\th'', t)\to X(\th'', t),
		\qquad \text{ as } k\to \infty. 
		\]
		For large enough $N$, we have
		\[
		m_- = \th_{i,N}< \th_{i+1,N}< \th' < m' <m''<\th'' < \th_{j-1,N}<\th_{j,N} = m_+.
		\] 
		Therefore 
		\[
		X_N(\th'',t)-X_N(\th',t)\leq X_N(\th_{j,N},t)-X_N(\th_{i+1,N},t)=x_{j,N}(t)-x_{i+1,N}(t)\leq D^0 e^{-\underline{\phi}t}.
		\]
		Taking $N\to \infty$ along the subsequence $(N_k)_{k=1}^\infty$ completes the proof.
	\end{proof}
	
	\subsection{Heavy-tailed and weakly singular communication: finite-time clustering}
	
	\label{ssec:wksing}
	
	As our final order of business, we consider \textit{weakly singular} communication protocols $\phi$, i.e., those satisfying \eqref{eq:phiws}.  These protocols are locally integrable but unbounded near the origin.  
		
	In contrast to the case of bounded communication $\phi$ (c.f. Theorem \ref{thm:critbound1}), finite-time cluster formation can occur in the critical region $\Sigma_0$ as a consequence of the weak singularity. This phenomenon was first discovered and analyzed in \cite{tan2020euler} for initial data satisfying $\Sigma_-=\emptyset$. The result is generalized in the following.
	
	\begin{THEOREM}\label{thm:critws}
		Assume $\phi$ is weakly singular (satisfying \eqref{eq:phiws}) and heavy-tailed (in particular, satisfying \eqref{eq:phimin}), and suppose that $L(m') = L(m'')$.  Then there exists a finite time $T$ such that $L(m')$ is a $t$-cluster for all $t\ge T$, i.e., $X(\cdot, t)$ is constant on $L(m')$ for $t\ge T$.  
	\end{THEOREM}
	
	\begin{proof}
		Denote $m_- = \inf L(m')$, $m_+ = \sup L(m')$, and assume without loss of generality that $m_- < m' < m'' < m_+$.  We prove that there exists a time $T$, which depends on $m_-$ and $m_+$ but not on $m'$, $m''$, such that $X(m',T) = X(m'',T)$.  This is enough to establish the theorem.  
		
		Choose a time $T_1$ such that $D^0e^{-\underline{\phi}T_1}=R$. We go through the proof of Theorem \ref{thm:critbound2} and deduce from \eqref{eq:xdiff5} that 
		\[x_{j,N}(t)-x_{i+1,N}(t)\leq R,\quad\forall~t\geq T_1.\]
		This allows us to apply \eqref{eq:phiws} and improve the estimate in \eqref{eq:phiconvrho5} by
		\begin{align*}
		&\sum_{\ell=1}^N m_{\ell,N}\int_{x_{i+1,N}(t)}^{x_{j,N}(t)}\phi(y-x_{\ell,N}(t))\dy\geq\sum_{\ell=i+1}^j m_{\ell,N}\int_{x_{i+1,N}(t)}^{x_{j,N}(t)}\phi(y-x_{\ell,N}(t))\dy\\
		&\qquad\geq \sum_{\ell=1}^N m_{\ell,N}\cdot c\big(x_{j,N}(t)-x_{i+1,N}(t)\big)^{-\b}\cdot\big(x_{j,N}(t)-x_{i+1,N}(t)\big)\\
		&\qquad=c(m_+-m_-)\big(x_{j,N}(t)-x_{i+1,N}(t)\big)^{1-\b}.	
		\end{align*}
		Hence, we have
		\[\frac{\dd}{\dt}\big(x_{j,N}(t)-x_{i+1,N}(t)\big)\leq 
		- c(m_+ - m_-) \big(x_{j,N}(t)-x_{i+1,N}(t)\big)^{1-\b},
		\qquad \forall t\ge T_1,
		\]
		and therefore
		\[x_{j,N}(t)-x_{i+1,N}(t)\leq \big(R^\b-c\b(m_+-m_-)(t-T_1)\big)^{1/\b},\quad\forall~t\geq T_1.\]
		It follows that $x_{j,N}(t)-x_{i+1,N}(t)$ becomes zero no later than
		\[T=T_1+\frac{R^\b}{c\b(m_+ - m_-)},\]
		where $T$ is independent of $N$.
		
		Finally, we proceed with the same argument as Step 2 in the proof of Theorem \ref{thm:critbound2} to pass to the limit and finish the proof.
	\end{proof}

\medskip
{\bf Acknowledgments.}
This material is based upon work supported by the National Science Foundation under Grant No. DMS-1928930 while TL was in residence at the Simons Laufer Mathematical Sciences Institute (formerly MSRI) in Berkeley, California, during the summer of 2023.

CT acknowledges the support of NSF grants DMS-2108264 and DMS-2238219.


\end{document}